\documentclass[12pt,a4paper]{amsart}

\usepackage[colorlinks=true,
linkcolor=green,
filecolor=brown,
citecolor=green]{hyperref}

\usepackage[
centering]{geometry}

\usepackage{amssymb}
\usepackage[euler-digits]{eulervm}
\usepackage{graphicx}
\usepackage[latin1]{inputenc}

\usepackage{parskip}

\title{A Quiver Presentation for Solomon's Descent Algebra.}
\author{Götz Pfeiffer}
\date{\today}
\address{Department of Mathematics, National University of Ireland, Galway,
University Road, Galway, Ireland}
\email{goetz.pfeiffer@nuigalway.ie}
\dedicatory{Dedicated to the memory of Manfred Schocker (1970--2006)}

\newcommand{\seqnum}[1]{\href{http://www.research.att.com/cgi-bin/access.cgi/as/~njas/sequences/eisA.cgi?Anum=#1}{\underline{#1}}}

\swapnumbers
\newtheorem{Theorem}[subsection]{Theorem}
\newtheorem{Proposition}[subsection]{Proposition}
\newtheorem{Lemma}[subsection]{Lemma}
\newtheorem{Corollary}[subsection]{Corollary}
\newtheorem{Conjecture}[subsection]{Conjecture}
\theoremstyle{definition}

\newtheorem{Example}[subsection]{Example}
\newtheorem{Remark}[subsection]{Remark}

\renewcommand{\labelenumi}{(\roman{enumi})}

\newcommand{\N}{\mathbb{N}}
\newcommand{\Q}{\mathbb{Q}}

\let\AA\relax\newcommand{\AA}{\mathcal{A}}
\newcommand{\BB}{\mathcal{B}}
\newcommand{\CC}{\mathcal{C}}
\newcommand{\DD}{\mathcal{D}}
\newcommand{\EE}{\mathcal{E}}
\newcommand{\GG}{\mathcal{G}}

\newcommand{\PP}{\mathcal{P}}

\newcommand{\RR}{\mathcal{R}}

\newcommand{\Qb}{\mathbf{Q}}
\newcommand{\Eb}{\mathbf{E}}
\newcommand{\eb}{\mathbf{e}}
\newcommand{\Vb}{\mathbf{V}}
\newcommand{\vb}{\mathbf{v}}

\DeclareMathOperator{\End}{\mathsf{End}}

\DeclareMathOperator{\Sym}{\mathsf{Sym}}

\DeclareMathOperator{\Rad}{\mathsf{Rad}}

\DeclareMathOperator{\id}{\mathsf{id}}

\DeclareMathOperator{\dd}{\mathsf{dp}}

\let\ker\relax\DeclareMathOperator{\ker}{\mathsf{ker}}
\let\dim\relax\DeclareMathOperator{\dim}{\mathsf{dim}}

\renewcommand{\emptyset}{\varnothing}
\let\coprod\bigsqcup

\renewcommand{\setminus}{-}

\newcommand{\prefix}{\preceq_{\pi}}
\newcommand{\suffix}{\succeq_{\sigma}}

\newcommand{\Size}[1]{\left|#1\right|}
\newcommand{\Span}[1]{\left<#1\right>}
\newcommand{\Floor}[1]{\lfloor#1\rfloor}

\newcommand{\inRCosets}[3]{#1 \in #3{/}#2}
\newcommand{\edge}[3]{#1\mathrel{\hbox to 0pt{\phantom{--}$^{#2}$\hss}{\longrightarrow}}#3}

\numberwithin{equation}{section}

\begin{document}

\begin{abstract}
  The descent algebra $\Sigma(W)$ is a subalgebra of the group algebra $\Q W$
  of a finite Coxeter group $W$, which supports a homomorphism with nilpotent
  kernel and commutative image in the character ring of $W$.  Thus
  $\Sigma(W)$ is a basic algebra, and as such it has a presentation as a
  quiver with relations.  Here we construct $\Sigma(W)$ as a quotient of a
  subalgebra of the path algebra of the Hasse diagram of the Boolean lattice
  of all subsets of $S$, the set of simple reflections in $W$.  From this
  construction we obtain some general information about the quiver of
  $\Sigma(W)$ and an algorithm for the construction of a quiver presentation
  for the descent algebra $\Sigma(W)$ of any given finite Coxeter group~$W$.
\end{abstract}

\maketitle

\section*{Contents.}
\setcounter{tocdepth}{1}
\tableofcontents

\section{Introduction.}

The descent algebra $\Sigma(W)$ of a finite Coxeter group $W$ of rank $n$ is
a remarkable $2^n$-dimensional subalgebra of the group algebra $\Q W$, which
supports a homomorphism $\theta$ with nilpotent kernel and commutative image
in the character ring of $W$.  Therefore, $\Sigma(W)$ is a basic algebra, and
as such it has a presentation as a quiver with relations.  In this article we
construct $\Sigma(W)$ as a quotient of a subalgebra of the path algebra of
the Hasse diagram of the power set of $S$, the set of simple reflections of
$W$.  From this construction we obtain general information about the quiver
of $\Sigma(W)$ and an algorithm, which for a given finite Coxeter group $W$
computes a quiver presentation for the descent algebra $\Sigma(W)$.  

Solomon~\cite{Solomon1976} has defined the descent algebra $\Sigma(W)$ in
terms of distinguished coset representatives of the standard parabolic
subgroups of $W$.  Bergeron, Bergeron, Howlett and Taylor~\cite{BeBeHoTa92}
have obtained a decomposition of $\Sigma(W)$ into principal indecomposable
modules.  More recently, Bidigare~\cite{bidigare1997} has identified
$\Sigma(W)$ with the fixed point space under the action of $W$ on the monoid
algebra of the face monoid of the hyperplane arrangement associated to the
reflection representation of $W$.  In this approach, the descent algebra
$\Sigma(W)$ is a subalgebra of a quotient of a path algebra.
Brown~\cite{brown2004} discusses this construction in the wider context of
semigroups of idempotents.  The face monoid algebra of the reflection
arrangement of a finite Coxeter group $W$ is called the Solomon-Tits algebra
by Patras and Schocker~\cite{PatrasSchocker2006}.
Schocker~\cite{Schocker2006} discusses the descent algebra of the symmetric
group and its quiver in this context.

More results concerning the quiver of $\Sigma(W)$ have been obtained for
particular types of finite Coxeter groups, mostly for type $A_n$.  Garsia and
Reutenauer~\cite{GarReut1989} have performed a very detailed analysis of the
descent algebra of the symmetric group, and described its quiver in terms of
restricted partition refinement.  Atkinson~\cite{Atkinson1992} has determined
the Loewy length of the descent algebra of the symmetric group.  Bonnaf\'e
and Pfeiffer~\cite{BonnafePfeiffer2008} have determined the Loewy length of
$\Sigma(W)$ for the other types of irreducible finite Coxeter groups with the
exception of type $D_n$ for  $n$ odd.
An argument which settles this case has been put forward by 
Saliola~\cite{saliola2007b}, based on his investigation of Bidigare's geometric setting~\cite{saliola2008}.

In this article, we present an alternative approach to $\Sigma(W)$.  We
construct a quiver with relations for $\Sigma(W)$ in three steps.  The point
of departure is the path algebra $A$ of the Hasse diagram of the power set
$\PP(S)$ of a finite set $S$, partially ordered by reverse inclusion.  Then
we use a partial action of $W$ on $\PP(S)$ to exhibit a subalgebra of $A$,
and a quiver presentation for it.  Finally, a quotient of this subalgebra,
formed with the help of a difference operator on $A$, is shown to be
isomorphic to $\Sigma(W)$.

This article is organized as follows.  In Section~\ref{sec:descents} we
recall the definition of the descent algebra in terms of the distinguished
coset representatives of parabolic subgroups of $W$ and some combinatorial
properties of these transversals.  Section~\ref{sec:quiver} introduces
quivers and their path algebras, and shows how monoid actions, in particular
of a free monoid, produce examples of quivers.  In Section~\ref{sec:shapes},
the conjugation action of $W$ on its parabolic subgroups is described as an
action of a free monoid on the standard parabolic subgroups of $W$.  In
Section~\ref{sec:alleys} we obtain the Hasse diagram of the power set of $S$
from the take-away action of the free monoid $S^*$.  The paths in this
particular quiver are called alleys and they form a basis of a path algebra
$A$.  Prefixes and suffixes of paths define in a natural way two rooted
forests on the set of all alleys.  In Section~\ref{sec:streets}, we apply the
conjugation action from Section~\ref{sec:shapes} to the alleys of
Section~\ref{sec:alleys}.  An orbit of alleys is called a street and the
streets (identified with the sums of their elements) form a basis of a
subalgebra $\Xi$ of $A$.  Prefixes and suffixes of alleys induce two rooted
forests on the set of all streets, which in particular decompose $\Xi$ into
projective indecomposable modules.  We furthermore conjecture that $\Xi$ is a
path algebra.  In Section~\ref{sec:difference}, a difference operator
$\Delta$ on $A$ is used to map $\Xi$ surjectively onto the grade $0$
component $A_0$ of $A$.  In Section~\ref{sec:matrix}, we use the difference
operator to define a matrix representation of $\Xi$ on $A_0$.  In
Section~\ref{sec:more-descents}, we prove in Theorem~\ref{thm:faxm} a key
result about right multiplication in $\Sigma(W)$.  We then identify $A_0$
with the descent algebra $\Sigma(W)$, and show as our Main
Theorem~\ref{thm:main} that with this identification $\Delta$ becomes an
anti-homomorphism from $\Xi$ onto $\Sigma(W)$.  In
Section~\ref{sec:properties}, we derive some properties of the quiver of
$\Sigma(W)$ from this construction.  Finally, in Section~\ref{sec:examples},
we present an algorithm which computes, for a given finite Coxeter group $W$,
a quiver presentation for $\Sigma(W)$.  For each of the series $A_n$, $B_n$,
and $D_n$ of irreducible finite Coxeter groups, we state some general
properties of the quiver of $\Sigma(W)$ and give one example of a quiver
presentation.

Throughout, we use the symmetric group on $4$ points for the purpose of
illustration.  The constructions, however, work for all types of finite
Coxeter groups.  Concrete results for particular types will be the subject of
subsequent articles.  Computer implementations of data structures
corresponding to some of the combinatorial and algebraic objects introduced
here have helped us to produce the examples and figures, and to verify
conjectured theorems in many cases.  They are available in the form of the
\textsf{GAP}~\cite{GAP} package \textsf{ZigZag}~\cite{zigzag}, which is based
on the \textsf{CHEVIE}~\cite{chevie} package for finite Coxeter groups and
Iwahori--Hecke algebras.

\section{Descents and Parabolic Transversals.} 
\label{sec:descents}

In this section, notation and some basic concepts are introduced, mostly
following Geck and Pfeiffer~\cite{GePf2000}.  Let $W$ be a finite Coxeter
group, generated by a set $S$ of simple reflections.  Let $\ell \colon W \to
\N_0$ be the usual length function on $W$.  The (left) \emph{descent set} of
an element $w\in W$ is the set
\begin{align}
\DD(w) = \{s \in S : \ell(sw) < \ell(w)\}.
\end{align}
For each subset
$J \subseteq S$, the subgroup $W_J = \Span{J}$ is called a (standard) 
\emph{parabolic subgroup} of $W$, and the set
\begin{align}
X_J = \{w \in W: \DD(w) \cap J = \emptyset\}
\end{align}
is a transversal of the right cosets $W_J w$ of the parabolic subgroup $W_J$
in $W$, consisting of the elements of minimal length in each coset.  For a
fixed subset $J \subseteq S$, each element $w \in W$ can be written as a
product $w = u \cdot x$ for unique elements $u \in W_J$ and $x \in X_J$.  A
product $w_1 w_2 \dotsm w_k$ of elements $w_1, w_2, \dots, w_k \in W$ is
called \emph{reduced} if
\begin{align} \label{eq:W-reduced}
  \ell(w_1 w_2 \dotsm w_k) = \ell(w_1) + \ell(w_2) + \dots + \ell(w_k).
\end{align}
If this is the case, we sometimes write a product like $w_1w_2w_3$ as $w_1
\cdot w_2 \cdot w_3$ in order to emphasize the fact that the product is
reduced.  For example, the product $u \cdot x$ of an element $u \in W_J$ and
a coset representative $x \in X_J$ is reduced.  The \emph{longest element} of
the parabolic subgroup $W_J$ is denoted by $w_J$, the longest element of $W$
also by $w_0$.  The elements $w_J$ are involutions.  The quotient $w_J^{-1}
w_0 = w_J w_0$ is the unique longest element of the transversal $X_J$.

The \emph{descent algebra} of $W$ is defined as the subspace $\Sigma(W)$ of
the group algebra $\Q W$ spanned by the sums
\begin{align}
x_J = \sum_{x \in X_J^{-1}} x \in \Q W
\end{align}
over the sets $X_J^{-1} = \{x^{-1} : x \in X_J\}$ 
of  \emph{left} coset representatives of $W_J$ in $W$,
for $J \subseteq S$.   By Solomon's Theorem~\cite{Solomon1976}, this subspace
is  in  fact  a  subalgebra  of   $\Q  W$  with  structure  constants  as  in
equation (\ref{eq:aJKL}) below.  For $J, K \subseteq S$, one further defines
\begin{align}
X_{JK} = X_J \cap X_K^{-1}
\quad \text{and} \quad
X_J^K = X_J \cap W_K.
\end{align}
The set $X_{JK}$ is a transversal of the double cosets of $W_J$ and $W_K$ in
$W$.  A parabolic transversal $X_J$ can in many ways be described in terms of
other transversals, or as a set of prefixes.  Here, an element $u \in W$ is
called a \emph{prefix} of $w \in W$ if $l(w) = l(u) + l(u^{-1} w)$.  In that
case we write $u \leq w$.  The partial order defined in this way on $W$ is
sometimes called the \emph{weak Bruhat order} on $W$.

\begin{Proposition}\label{pro:xJ}
  Let $J, K \subseteq S$. Then
  \begin{enumerate}
  \item $X_J = X_J^K \cdot X_K$ if $J \subseteq K$;
  \item $X_J = d X_K$ if $d \in X_J$ and $K$ are such that $J^d = K$;
  \item $X_J = \coprod_{d \in X_{JK}} d \cdot X_{J^d \cap K}^K$;
  \item $X_J = \{w \in W : w \leq w_J w_0\}$.
Thus, $w \in X_J$ whenever $w \leq x$ for some $x \in X_J$.
  \end{enumerate}
\end{Proposition}

\begin{proof}
  \cite[(2.1.5), (2.1.8), (2.1.9), and (2.2.1)]{GePf2000}.
\end{proof}

For subsets $J, K, L \subseteq S$, we furthermore define
\begin{align} \label{eq:xJKL}
X_{JKL} = \{x \in X_{JK} : J^x \cap K = L\}.
\end{align}
The cardinalities 
\begin{align}
  a_{JKL} = \Size{X_{JKL}}
\end{align}
are the structure constants of the descent algebra: according to
Solomon~\cite{Solomon1976}, for all $J, K \subseteq S$,
\begin{align} \label{eq:aJKL}
x_J\, x_K = \sum_{d \in X_{JK}} x_{J^d \cap K} = \sum_{L \subseteq S} a_{JKL}
\, x_L.
\end{align}

Denote by $\theta$ the linear map from $\Sigma(W)$ into the character ring of
$W$ (over $\Q$) which assigns, for $J \subseteq S$, to $x_J$ the permutation
character 
\begin{align}
  \theta(x_J) = 1_{W_J}^W
\end{align}
of the action of $W$ on the cosets of the
parabolic subgroup $W_J$ in $W$.  Then, according to Solomon~\cite{Solomon1976},
$\theta$ is a homomorphism of algebras with commutative image and nilpotent
kernel.  It follows that the descent algebra $\Sigma(W)$ is a basic algebra,
and as such it has a presentation as a quiver with relations.

\section{Quivers, Path Algebras and Monoid Actions.} 
\label{sec:quiver}

A \emph{quiver} is a directed multigraph $Q = (V, E)$ consisting of a vertex
set $V$ and an edge set $E$, together with two maps $\iota, \tau \colon E \to
V$, assigning to each edge $e \in E$ a \emph{source} (or \emph{initial vertex})
$\iota(e) \in V$ and a \emph{target} (or \emph{terminal vertex}) $\tau(e) \in V$.
A \emph{path} of \emph{length} $\ell(a) = l$ in $Q$ is a pair 
\begin{align} \label{eq:path}
  a = (v; e_1, e_2, \dots, e_l)
\end{align}
consisting of a source $v \in V$ and a sequence of $l$ edges
$e_1, e_2, \dots, e_l \in E$ such that $\iota(e_1) = v$ and $\iota(e_i) =
\tau(e_{i-1})$ for $i = 2, \dots, l$.

Let $\AA$ be the set of all paths in $Q$ and let 
\begin{align}
  \AA_l = \{a \in \AA: \ell(a) = l\}
\end{align}
be the set of paths of length $l$.  We denote by $\emptyset$
the unique sequence of length $0$ and identify a
vertex $v \in V$ with the path $(v; \emptyset) \in \AA_0$.  We also identify
an edge $e \in E$ with the path $(\iota(e); e) \in \AA_1$.  The following
properties of path concatenation are obvious.

\begin{Proposition} \label{prop:a-category} The set $\AA = \coprod_{l \geq 0}
  \AA_l$ together with the partial multiplication defined as
\begin{align*}
  (v; e_1, \dots, e_l) \circ (v'; e_1', \dots, e_{l'}')
= (v; e_1, \dots, e_l, e_1', \dots, e_{l'}'),
\end{align*}
provided that $\tau(e_l) = v'$, is (the set of morphisms of) a category with
object set $\AA_0$.  The category $\AA$ is the free category generated by the
quiver $Q$.  Every path $a \in \AA$ of length $\ell(a) > 0$ is a unique
product of elements in the set $\AA_1$.
\end{Proposition}

The category $\AA$ of all paths in $Q$ can be used as formal basis of a
vector space.  For $l \geq 0$, let $A_l = \Q[\AA_l]$, the $\Q$-vector space
with basis $\AA_l$.  The \emph{path algebra} $A$ of the quiver $Q$ is defined
as
\begin{align}
  A = \Q[\AA] = \bigoplus_{l\geq 0} A_l,
\end{align}
where $a \circ a' = 0$ if the product $a \circ a'$ is not defined in $\AA$,
and otherwise multiplication is extended by linearity from $\AA$.  The path
algebra $A$ is a graded algebra, since we have $A_l \circ A_k \subseteq
A_{l + k}$, for all $l, k \geq 0$.

\subsection{Monoid Actions.} \label{sec:monoid}

A good source for examples of quivers and categories are monoid actions.
Suppose $M$ is a monoid acting on a set $X$ via $(x, m) \mapsto x.m$, then
the set $X \times M$ together with the partial multiplication
\begin{align}
(x,  m) \circ (x', m') = (x, mm'),
\end{align}
whenever $x, x' \in X$ and $m, m' \in M$ are such that $x.m = m'$, is a
category with object set $X$.

If $M$ is generated by a set $S \subseteq M$, then the \emph{action graph}
defined as the directed multigraph with vertex set $X$ and edge set $X \times
S$ is a quiver with $\iota(x, s) = x$ and $\tau(x, s) = x.s$ for all $(x, s)
\in X \times S$.  If $M$ is the free monoid $S^*$, then the category $X
\times M$ is the free category generated by the action graph $(X, X \times
S)$.

In the following two sections, we consider two different, but related,
examples of actions of the free monoid $S^*$ on the power set $\PP(S)$ of a
finite set~$S$.  In section~\ref{sec:streets}, we apply one action to the
path algebra arising from the other.

\section{Shapes.} \label{sec:shapes}

Let $W$ be the finite Coxeter group from Section~\ref{sec:descents} and let
$S \subseteq W$ be its set of simple reflections.  Here we regard $W$ as a
quotient of the free monoid $S^*$ consisting of all words over the alphabet
$S$.  The empty word will be denoted by $\emptyset$; the identity element of
$W$ by $\id$.  Beyond that, we make no notational effort to distinguish a
word in $S^*$ from a product of simple reflections in $W$.

The conjugation action  of $W$ on itself induces a  conjugation action of $W$
on its subsets which partitions  the power set $\PP(S) \subseteq \PP(W)$ into
classes of $W$-conjugate subsets.  We write  
\begin{align}
  A \sim B
\end{align}
if $A, B \subseteq W$
are such that $B = A^w$ for some $w \in W$ and call the class
\begin{align}
  [J] = \{K \subseteq S : K \sim J\}
\end{align}
of a subset $J \subseteq S$ the \emph{shape} of $J$ in $W$.  Moreover, we
denote by
\begin{align}
  \Lambda = \{[J] : J \subseteq S\}
\end{align}
the set of all shapes of $W$.  The shapes of $W$ correspond to the conjugacy
classes of parabolic subgroups of~$W$, since, for $J, K \subseteq S$, the
parabolic subgroup $W_K$ is a conjugate of $W_J$ if and only if $K$ is a
conjugate of $J$~\cite[(2.1.13)]{GePf2000}.

\begin{Example}
  Suppose $W$ is the symmetric group $\Sym(n+1)$ of degree $n+1$.  Every
  parabolic subgroup of $W$ is a direct product of symmetric groups, whose
  degrees form a partition of $n+1$, if fixed points are counted as factors
  of degree~$1$.  Two parabolic subgroups are conjugate in $W$ if and only if
  the corresponding partitions are the same.  In this way the shapes of a
  Coxeter group of type $A_n$ correspond to the partitions of $n+1$.

  In a similar way, the shapes of a Coxeter group of type $B_n$ correspond to
  the partitions of all $m \in \{0, \dots, n\}$~\cite[(2.3.10)]{GePf2000}.
  The shapes of a Coxeter group of type $D_n$ correspond to the partitions of
  all $m \in \{0, \dots, n-2\}$, and the partitions of $n$, with two copies
  of each even partition of $n$~\cite[(2.3.13)]{GePf2000}.
\end{Example}

In order to decide whether a parabolic subgroup $W_K$ of $W$ is a conjugate
$W_J^x$ of a parabolic subgroup $W_J$, it clearly suffices to consider
elements $x \in X_J$ which conjugate $J \subseteq S$ to a subset of $S$.  For $J
\subseteq S$ we denote
\begin{align} \label{eq:xjhash}
X_J^{\sharp} = \{x \in X_J : J^x \subseteq S\} = \coprod_{K \sim J} X_{JKK}.
\end{align}
Certainly $\id \in X_J^{\sharp}$ for all $J \subseteq S$.  And it is easy to
see that $X_L^{\sharp} \subseteq X_J^{\sharp}$ for all $J \subseteq L
\subseteq S$.  In general, the elements of $X_J^{\sharp}$ can be described as
reduced products of certain longest coset representatives.  Given $J
\subseteq S$ and $s \in S$, we let $L = J \cup \{s\}$ and  denote
\begin{align} \label{eq:omega}
  \omega(J, s) = w_J w_L \in X_J^{\sharp},
\end{align}
the longest coset representative of $W_J$ in $W_L$.  Clearly, if $s \in J$
then $\omega(J, s) = \id$.  Note that
\begin{align} \label{eq:omega-inverse}
  \omega(J, s)^{-1} = \omega(J^{w_L}, s^{w_L}).
\end{align}

\begin{Lemma} \label{la:sharp} 
  Let $J \subseteq S$ and let $x \in X_J^{\sharp}$.
  \begin{enumerate}
  \item $x^{-1} \in X_{J^x}^{\sharp}$.
  \item  If $J$ is a maximal subset of $S$ then $x = \id$ or $x = w_J w_0$.
  \item If $s \in \DD(x)$ then $\omega(J, s)$ is a prefix of~$x$.  Moreover, ${\omega(J, s)^{-1} x \in X_K^{\sharp}}$, where $K =
    J^{\omega(J, s)}$.
  \item There exist elements $s_1, \dots, s_r \in S$ such that
    \begin{align*}
      x = \omega(J_1, s_1) \cdot \omega(J_2, s_2) \dotsm \omega(J_r, s_r)
    \end{align*}
is a reduced product, where  $J_1 =
J$ and $J_{k+1} = J_k^{\omega(J_k, s_k)}$ for $1 \leq k < r$.
\end{enumerate}
\end{Lemma}

\begin{proof}
  (i) If $J^x = K \subseteq S$ then clearly $K^{x^{-1}} = J \subseteq S$.
  It remains to show that $x^{-1} \in X_K = \{w \in W : K \cap \DD(w) =
  \emptyset\}$.  But if there exists an element $s \in K \cap
  \DD(x^{-1})$  then $s^{x^{-1}} = xsx^{-1} \in J \cap \DD(x)$,
  contradicting $x \in X_J$.

  (ii) \cite[(2.3.2)]{GePf2000}. 

  (iii) Let $L = J \cup \{s\}$.  Using Proposition~\ref{pro:xJ}(i), we can
  write $x = x_1 \cdot x_2$ for (unique) elements $x_1 \in X_J^L$ and $x_2
  \in X_L$.  Clearly $x_1 \in X_J^{\sharp}$ and $x_1 \neq \id$ since $s \in
  \DD(x_1)$.  With (ii), this forces $x_1 = w_Jw_L$, since $J$ is a maximal
  subset of $L$.

  (iv) follows by induction on $\ell(x)$ from (iii).
\end{proof}

The preceding lemma motivates the following definition of 
an action of the
free monoid $S^*$ on the power set $\PP(S)$.
For $J \subseteq S$
and $s \in S$   we set
\begin{align}
  J.s = J^{\omega(J, s)}.
\end{align}

Most of the following properties are obvious.

\begin{Lemma} \label{la:J.s} 
Let $J \subseteq S$, $s \in S$ and $L = J \cup \{s\}$.  Then
\begin{enumerate}
\item $J.s = J$ if $s \in J$;
\item $J.ss = J.s$;
\item $s^{w_L} \in S$ and $J.s s^{w_L} = J$.
\end{enumerate}
\end{Lemma}

\begin{proof}
  (ii) If $J.s \neq J$ then $s
  \in J.s$.  (iii) See equation~(\ref{eq:omega-inverse}).
\end{proof}

Lemma~\ref{la:J.s}(iii) shows that all the effects of $S^*$ on $\PP(S)$ can
be undone.  Therefore, the $S^*$-orbits on $\PP(S)$ form a partition of
$\PP(S)$.  In fact, as a consequence of Lemma~\ref{la:sharp}(iv), these orbits
coincide with the shapes of $W$.  The following theorem states some known
properties of normalizers of parabolic subgroups in terms of the conjugation
action of $S^*$ on the power set $\PP(S)$.

\begin{Theorem}[Howlett]
  Suppose $S^*$ acts on $\PP(S)$ as defined above.
  \begin{enumerate}
  \item The orbits of $S^*$ on $\PP(S)$ form a partition of $\PP(S)$ and, for
    all $J \subseteq S$, the $S^*$-orbit of $J$ coincides with the shape
    $[J]$ in~$W$.
  \item The stabilizer 
    \begin{align*}
      N_J = \{x \in X_J : J^x = J\} = X_{JJJ}
    \end{align*}
    of $J$ is a complement of $W_J$ in its normalizer $N_W(W_J) = W_J \rtimes
    N_J$.
  \end{enumerate}
\end{Theorem}

\begin{proof}
  \cite{Howlett1980}; see also \cite[Theorem 2.3.3 and Proposition
  2.1.15]{GePf2000}.
\end{proof}

According to Section~\ref{sec:monoid}, this action gives rise to a category
$\PP(S) \times S^*$.  We will ignore the obviously trivial elements $(J,s)$
where $s \in J$ (see Lemma~\ref{la:J.s}(i)) and  define the
\emph{action graph} of $S^*$ on $\PP(S)$ as the quiver with edge set 
\begin{align}
  \{(J, s) : J \subseteq S,\, s \in S\setminus J\},
\end{align}
where $\iota(J, s) = J$ and $\tau(J, s) = J.s$.  We then define a category
$\CC$ as the subcategory of $\PP(S) \times S^*$ which is generated as a free
category by this action graph.

\begin{Example}
  Figure~\ref{fig:shape3} illustrates the action graph and the category $\CC$
  for the type~$A_3$ generated by a set $S = \{1,2,3\}$ of simple
  reflections, where $1$ commutes with $3$.  Here the vertices are the
  subsets of $S$, written with braces and commas omitted.  Multiple labels
  (like $1,2,3$) on arrows indicate multiple arrows, one for each label.
  \begin{figure}[htbp]
    \centering
    \includegraphics[scale=.3]{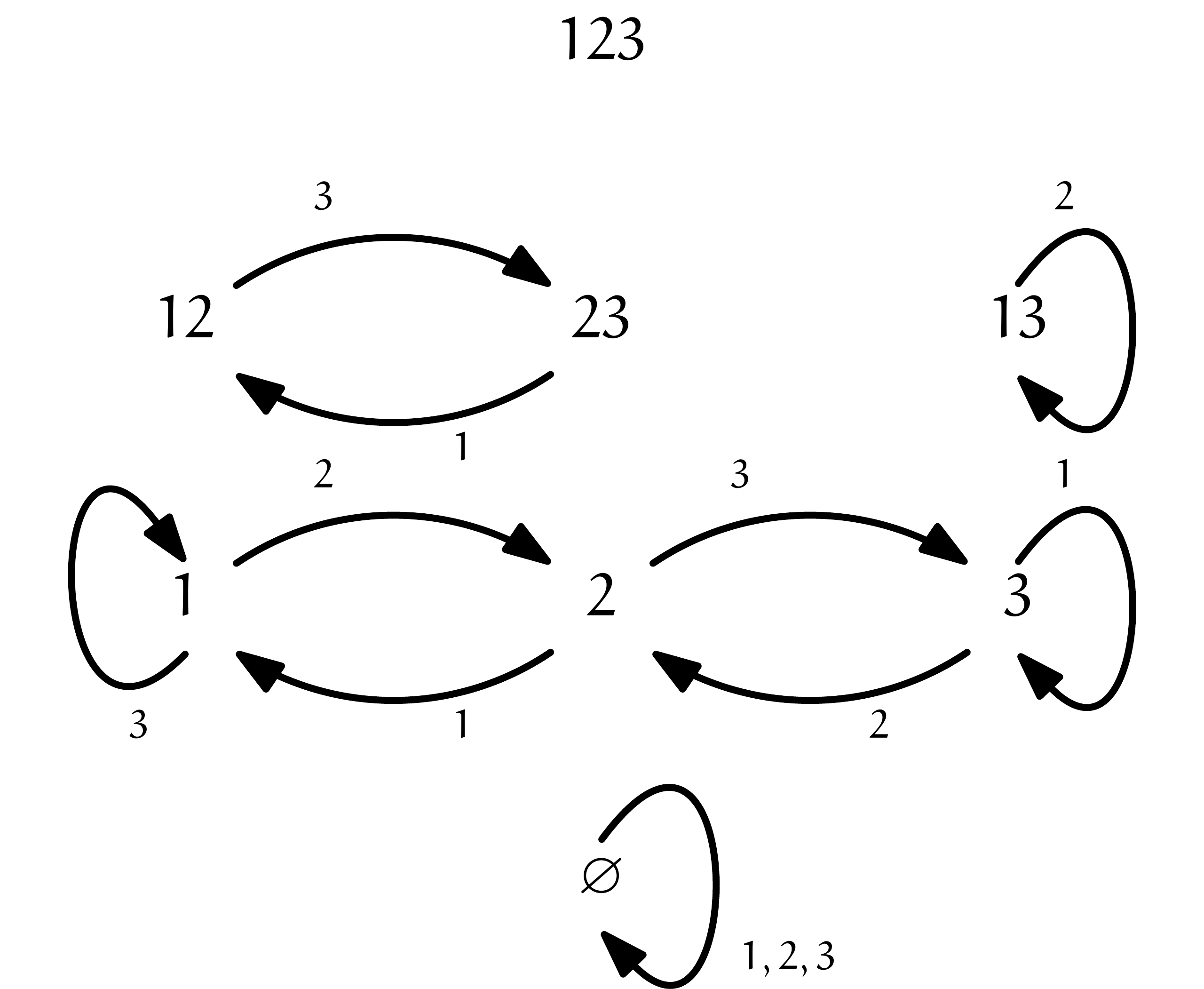}
    \caption{$S^*$ acting on $\PP(S)$ of type $A_3$.}
    \label{fig:shape3}
  \end{figure}
\end{Example}

The fact that the effects of the monoid $S^*$ on any subset $J \subseteq S$
can be undone manifests itself further in the form of a closely related
groupoid.  Note that if $J \sim K$ then $X_{JKK} = \{x \in X_J : J^x = K\}$.
From Proposition~\ref{pro:xJ}(ii) it follows that, if $J \sim K \sim L$ then
$d X_{KLL} = X_{JLL}$ for all $d \in X_{JKK}$, and hence that $X_{JKK}
X_{KLL} = X_{JLL}$.  Recall that $X_J^{\sharp} = \coprod_{K \sim J} X_{JKK}$.
Then the set of pairs
\begin{align}
 \GG = \{(J, x) : J \subseteq S,\, x \in X_J^{\sharp}\}
\end{align}
forms a category with respect to the partial multiplication 
\begin{align}
  (J, x) (K, y) = (J, xy)
\end{align}
if $J, K \subseteq S$, $x \in X_J^{\sharp}$ and $y \in X_K^{\sharp}$ are such
that $J^x = K$.  Each pair $(J, x) \in \GG$ has an inverse $(J^x, x^{-1}) \in
\GG$, since  by Lemma~\ref{la:sharp}(i) $x \in  X_J^{\sharp}$ implies $x^{-1}
\in X_{J^x}^{\sharp}$.  Therefore, the category $\GG$ is in fact a groupoid.

There exists a unique functor $\omega$ from the category $\PP(S) \times S^*$
to the group $W$, regarded as a one-object category, extending the map
$\omega \colon \PP(S) \times S \to W$ as defined in
equation~(\ref{eq:omega}).  Then, 
\begin{align}
  J.m = J^{\omega(J, m)} 
\end{align}
for all $J \in
\PP(S)$ and all $m \in S^*$, and conjugation by $\omega(J, m)$ induces a
bijection
\begin{align} \label{eq:partialbijection}
  s \mapsto s^{\omega(J, m)}
\end{align}
from $J$ to $J.m$.

\begin{Proposition} \label{pro:C-to-G}
  The map $(J, m) \mapsto (J, \omega(J, m))$ together with the identity map
  on the object set $\PP(S)$ is a functor from the category $\PP(S) \times
  S^*$ onto the groupoid~$\GG$.
\end{Proposition}

\begin{proof}
  Let $J \subseteq S$.  We only need to show that $\omega(J, m) \in X_J$ for
  all $m \in S^*$.  If $m = \emptyset$ then $\omega(J, m) = \id$ and there is
  nothing to show.  Otherwise $m = sm'$ for some $s \in S$ and $m' \in S^*$.
  Let $d = \omega(J, s)$.  Then $\omega(J, m) = d\omega(J^d, m')$, where
  $\omega(J^d, m') \in X_{J^d}$, by induction on the length of $m'$.  By
  Proposition~\ref{pro:xJ}, $d X_{J^d} = X_J$, whence $\omega(J, m) =
  d\omega(J^d, m') \in X_J$.
\end{proof}

Brink and Howlett~\cite{BrinkHowlett1999} have given a presentation of $\GG$
as a quotient of the free category $\CC$ in terms of the generating set
\begin{align}
\{(J, \omega(J, s)) : J \subseteq S, s \in S \setminus J\}
\end{align}
corresponding to the edges of the action graph.  Before we formulate the
relations, we introduce the notion of a reduced expression for elements of
$\GG$.

We call $(J, m) \in \CC$ with $m = s_1 s_2 \dotsm s_k \in S^*$ a
\emph{reduced expression} for $(K, x) \in \GG$ if $(J, \omega(J, m)) = (K, x)$
and
\begin{align}
  x = \omega(J_1, s_1) \cdot
\omega(J_2, s_2) \dotsm
\omega(J_k, s_k)
\end{align}
is a reduced product in $W$ in the sense of equation~(\ref{eq:W-reduced}),
where $J_i = J.s_1 s_2 \dotsm s_{i-1}$, for $i = 1, \dots, k$.  In
particular, for $J \subseteq S$, the pair $(J, \emptyset) \in \CC$ is the
unique reduced expression for $(J, \id) \in \GG$.

By Lemma~\ref{la:sharp}(iv), each element $(J, x) \in \GG$ has a reduced
expression.  Suppose $J \subseteq L \subseteq M \subseteq S$.  Then clearly
$w_L w_M \in X_J^{\sharp}$ and $(J, w_L w_M) \in \GG$ has a reduced
expression.  The following result by Brink and
Howlett~\cite{BrinkHowlett1999} is concerned with a situation where this
reduced expression is unique.

\begin{Proposition}
  Let $J \subseteq S$ and $s, t \in S \setminus J$ be such that $s \neq t$.
Let $L = J \cup \{s\}$ 
 and let $w \in X_J^{\sharp}$ be a prefix of
  $\omega(L, t)$.  Then the element $(J, w) \in \GG$ has a unique reduced
  expression.  In particular, the reduced expression for $(J, \omega(L, t)) \in
  \GG$ is unique.
\end{Proposition}

\begin{proof}
  Without loss of generality we may assume that $L \cup \{t\} = S$.  Then, by
  Proposition~\ref{pro:xJ}(iv), $w$ is a prefix of $\omega(L, t)$ if and only
  if $w \in X_L$.

  The claim is certainly true for $w = \id$.  Otherwise, $w \in X_J^{\sharp}
  \cap X_L$ is such that $l(w) > 0$.  If follows that $\DD(w) = \{t\}$ and
  hence, using Lemma~\ref{la:sharp}(iii), that each reduced expression for
  $(J, w) \in \GG$ begins with $(J, t) \in \CC$.

  Let $K = J \cup \{t\}$.  A straightforward comparison of lengths shows that
  the product
  \begin{align}
    (\omega(J, t)^{-1} \omega(L, t)) \cdot \omega(J, s)^{\omega(L, t)} &=
 ((w_J w_K)^{-1} w_L w_0) \cdot (w_J w_L)^{w_L w_0} 
\\\notag&= w_K w_0 = \omega(K, s)
  \end{align}
  is reduced in $W$.  Therefore the element
  \begin{align}
    w' = \omega(J, t)^{-1} w \in X_{J.t}^{\sharp}
  \end{align}
  is a prefix of $\omega(K, s)$.  By induction on the length $\ell(w)$, the pair
  $(J.t, w') \in \GG$ has a unique reduced expression.  This shows that
  \begin{align}
    (J, w) = (J, \omega(J, t)) \cdot  (J.t, w') \in \GG
  \end{align}
  has a unique reduced expression.
\end{proof}

The last result implies that $(J, w_J w_M)$, where $M = J \cup \{s, t\}$, has
exactly two reduced  expressions, $(J, s m_L)$ and $(J,  t m_K)$, where $m_L,
m_K  \in S^*$  are  such that  $(J.s, m_L)  \in  \CC$ is  the unique  reduced
expression for $(J.s, \omega(L, t)) \in  \GG$ and $(J.t, m_K) \in \CC$ is the
unique reduced expression for $(J.t, \omega(K, s)) \in \GG$.
These two reduced expressions for the same element give rise
to a relation
\begin{align} \label{eq:rel2}
  (J, s m_L) \equiv (J,  t m_K).
\end{align}
It is a remarkable fact, that these relations together
with relations of the form
\begin{align} \label{eq:rel1}
  (J, s s^{w_L}) \equiv (J, \emptyset),
\end{align}
which are consequences of (\ref{eq:omega-inverse}),
are sufficient for a presentation of  the groupoid $\GG$ as a quotient of the
free category $\CC$.

\begin{Theorem}[Brink--Howlett]
Let $\equiv$ be the congruence generated by all the relations of the
forms (\ref{eq:rel1}) and (\ref{eq:rel2}) in $\CC$.  Then
  $\CC/{\equiv}$ is isomorphic to $\GG$.
\end{Theorem}

\section{Alleys.} \label{sec:alleys}

In this section we introduce a particular quiver $Q$ on the vertex set
$\PP(S)$, whose paths will be called alleys.  The quiver $Q$ arises from the
take-away action of the free monoid $S^*$ on $\PP(S)$ defined by $(L, s)
\mapsto L \setminus \{s\}$.  We exhibit various structural properties of the
category $\AA$ of all alleys and of the corresponding path algebra $A$ of the
quiver $Q$.  There are two natural partial orders on $\AA$, which turn the
set $\AA$ into a rooted twofold forest.  In subsequent sections, the
$S^*$-action on $\PP(S)$ from Section~\ref{sec:shapes} will be extended to an
$S^*$-action on $\AA$ and a difference operator~$\delta$ will be defined on
the graded algebra~$A$.

Let $S$ be a finite set.  For $L \subseteq S$ and $s \in S$ denote
\begin{align}
  L_s = L \setminus \{s\}.
\end{align}
The map $\PP(S) \times S \to \PP(S),\, (L, s) \mapsto L_s$, defines an action
of the free monoid $S^*$ on the power set $\PP(S)$.  Ignoring the obviously
trivial edges $(L, s)$ where $s \in S \setminus L$, the action graph of
this action is the Hasse diagram $Q$ of the power set $\PP(S)$, partially
ordered by reverse inclusion.  For a subset $L \subseteq S$ and pairwise
different elements $s, t, \dotsc \in L$ denote by $(L; s, t, \dots)$ the
unique path with vertices
\begin{align}
  L,\,
 L \setminus \{s\},\,
 L \setminus \{s, t\},\, \dots,\,
 L \setminus \{s, t, \dots\}
\end{align}
in $Q$, deviating slightly from the notation for quivers in
equation~(\ref{eq:path}).  As before in Section~\ref{sec:quiver}, we denote
by $\AA$ the category of all paths in $Q$, and by $A = \Q[\AA]$ the path
algebra of the quiver $Q$.  We call the path $a = (L; s, t, \dots) \in \AA$
an \emph{alley} from its \emph{source} $\iota(a) = L$ to its \emph{target}
$\tau(a) = L \setminus \{s, t, \dots\}$.  The elements of the sequence $(s,
t, \dots)$ are the \emph{segments} of $a$, and the \emph{length} $\ell(a)$ of
$a$ is the number of segments, i.e., $\ell(a) = \#\{s, t, \dots\}$.

\begin{Example}
  Figure~\ref{fig:power3i} shows the quiver $Q$ for $S = \{1, 2, 3\}$.  As in
  Figure~\ref{fig:shape3}, the vertices of $Q$ are the subsets of $S$,
  written without punctuation.
  \begin{figure}[htb]
    \centering
    \includegraphics[scale=.3]{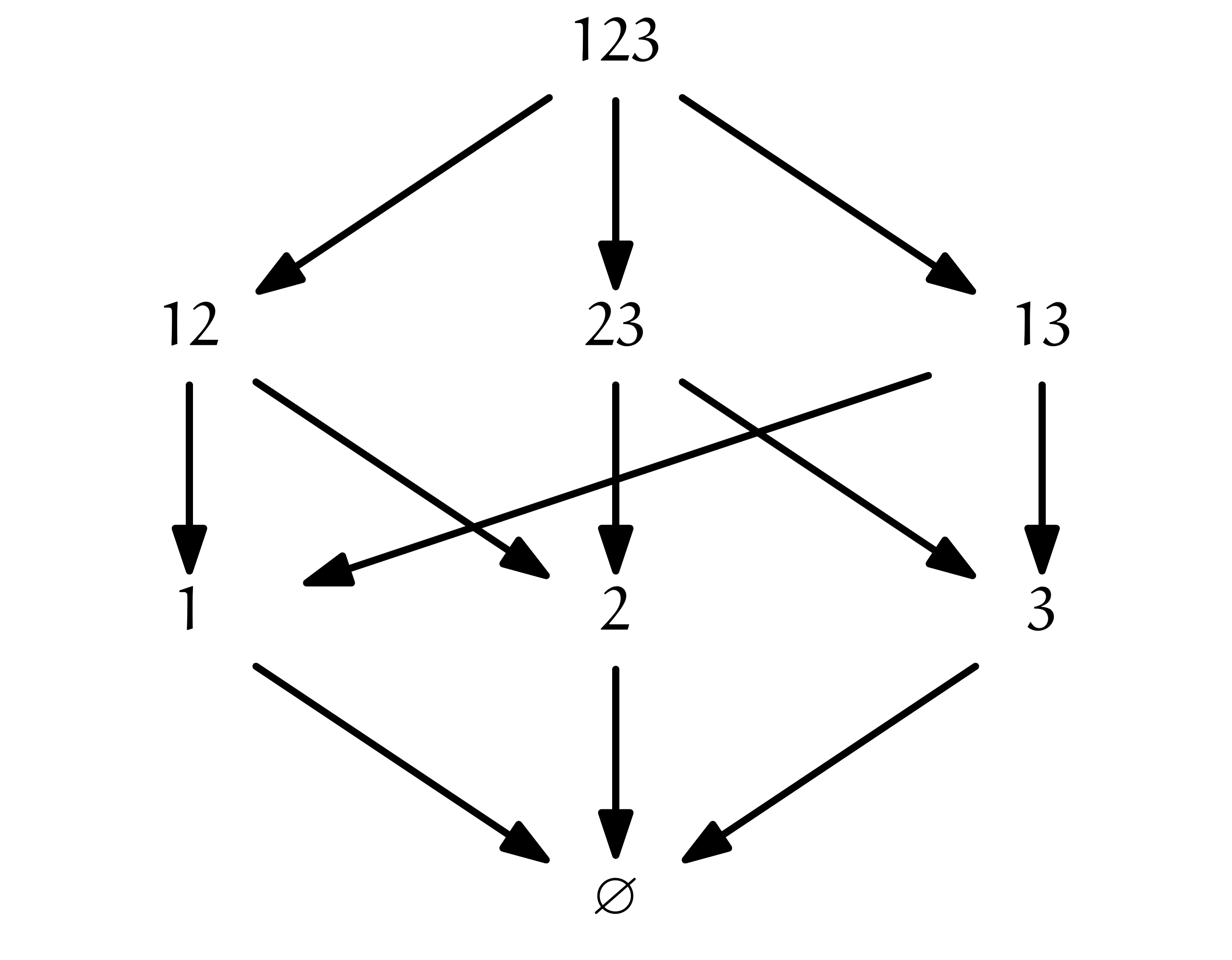}
    \caption{The Alley Quiver $Q$ for $S = \{1, 2, 3\}$.}
    \label{fig:power3i}
  \end{figure}
\end{Example}

\subsection{Counting Alleys.} 

Let $n = \Size{S}$.
For $0 \leq l \leq n$ there are $\frac{2^l}{l!}n!$ alleys of length $n-l$,
since there are $\frac{n!}{l!}$ ways to choose $n - l$ segments and $2^l$
ways to choose a target from the remaining $l$ elements of $S$.  Hence there
are
\begin{equation}
\Size{\AA} = n! \sum_{l=0}^n \frac{2^l}{l!}
\end{equation}
alleys in total.  The values of $\Size{\AA}$ for $n \leq 9$ are as follows.
\begin{align}
\begin{array}{c|cccccccccc}
  n & 0 & 1 & 2 & 3 & 4 & 5 & 6 & 7 & 8 & 9 \\ \hline
\Size{\AA} & 1 & 3 & 10 & 38 & 168 & 872 & 5296 & 37200 & 297856 & 2681216
\end{array}
\end{align}
This is  sequence number \seqnum{A010842} in Sloane's  online encyclopedia of
integer  sequences~\cite{OEIS},  which  has exponential  generating  function
$e^{2x}/(1-x)$.

\subsection{Partial Order.} \label{sec:partial-order}

The set $\AA$ of all alleys is in two ways partially ordered as follows.  Let
$a, a' \in \AA$.  We say that  $a'$ is a \emph{prefix} of $a$, and write
\begin{align}
a' \prefix a,
\end{align} 
if $a = a' \circ a''$ for some $a'' \in \AA$.  By the unique factorization
property of Proposition~\ref{prop:a-category}, each $a = (L; s_1, \dots, s_l)
\in \AA$ of length $\ell(a) = l > 0$ has a unique longest nontrivial prefix
\begin{align}
\pi(a) = (L; s_1, \dots, s_{l-1}).
\end{align}
Thus $\AA$ is a forest of rooted trees with roots $L \subseteq S$, where
$\pi(a)$ is the parent of $a$ if $\ell(a) > 0$.  The alley $a$ lies in the
$\pi$-tree with root $L$ if and only if $\iota(a) = L$.  The
$\pi$-\emph{children} of an alley $a = (L; s_1, s_2, \dots, s_l)$ are the
alleys $(L; s_1, s_2, \dots, s_{l+1})$ with $s_{l+1} \in L \setminus \{s_1,
\dots, s_l\}$.  The alley $a = (L; s_1, s_2, \dots, s_l)$ is a
$\pi$-\emph{leaf} if $L = \{s_1, \dots, s_l\}$.  Note that
\begin{align}
  \iota(a) = \pi^{\ell(a)}(a)
\end{align}
for all $a \in \AA$.

We furthermore  say that $a'$ is a  \emph{suffix} of $a$, or  that $a$ \emph{ends
  in} $a'$, and write
\begin{align}
a \suffix a'
\end{align}
if $a = a'' \circ a'$ for some $a'' \in \AA$. By the unique factorization
property of Proposition~\ref{prop:a-category}, each $a = (L; s, t, \dots)
\in \AA$ of length $\ell(a) > 0$ has a unique longest nontrivial suffix
\begin{align}
\sigma(a) = (L_s; t, \dots).
\end{align}
Thus $\AA$ also is a forest of rooted trees with roots $L \subseteq S$, where
$\sigma(a)$ is the parent of $a$ if $\ell(a) > 0$.  Here the alley $a$ lies
in the $\sigma$-tree with root $L$ if and only if $\tau(a) = L$.  The
$\sigma$-\emph{children} of an alley $a = (L; s, t, \dots)$ are the alleys $(L
\cup \{r\}; r, s, t, \dots)$ with $r \in S \setminus L$.  The alley $a = (L;
s, t, \dots)$ is a $\sigma$-\emph{leaf} if $L = S$. Note that
\begin{align}
  \tau(a) = \sigma^{\ell(a)}(a)
\end{align}
for all $a \in \AA$.

Hence $\AA$ together with the two partial orders $\prefix$ and $\suffix$
forms a rooted twofold forest where the subsets $L \subseteq S$ serve as roots
for both forests.  In terms of these forests, the product of two alleys $a,
a' \in \AA$ can be described as the intersection of the $\pi$-subtree with root
$a$ and the $\sigma$-subtree with root $a'$, provided that $\tau(a) =
\iota(a')$; in that case $\tau(a) = \iota(a')$ is the only subset $L
\subseteq S$ with $a \suffix L \prefix a'$, otherwise there is no such subset
at all.  In any case, this argument proves the following result.

\begin{Proposition} \label{pro:alley-forest-product}
  Let $a, a' \in \AA$.  Then
\begin{align*}
a \circ a' = \mathop{\sum\sum}_{\substack{L \subseteq S,a'' \in \AA \\a \suffix L \prefix a'\\a \prefix a'' \suffix a'}} a''.
\end{align*}
\end{Proposition}

\section{Streets.} \label{sec:streets}

As before in Section~\ref{sec:shapes}, let $S$ be the set of simple
reflections of a finite Coxeter group~$W$.  In this section we relate the
conjugation action of the free monoid $S^*$ on the power set $\PP(S)$ from
Section~\ref{sec:shapes} to the takeaway action of $S^*$ on $\PP(S)$ from
Section~\ref{sec:alleys}, by extending the conjugation action to the set
$\AA$ of all alleys.  Again, the $S^*$-orbits form a partition of $\AA$.  An
orbit of alleys will be called a street and $\Psi$ will denote the set of all
streets.  The two main results of this section show that $\Psi$ is a rooted
twofold forest, and that the linear span of $\Psi$ forms a subalgebra $\Xi$
of~$A$.  The partial order on streets allows us to identify a complete set of
primitive orthogonal idempotents for $\Xi$.  We furthermore conjecture that
$\Xi$ is a path algebra.

The conjugate of $a = (L; s, t, \dots) \in \AA$ by an
element $w \in W$ is the pair
\begin{align}
  a^w = (L^w; s^w, t^w, \dots),
\end{align}
consisting of a subset $L^w \subseteq W$ and a sequence of elements $s^w,
t^w, \dotsc \in L^w$.  The conjugate $c^w$ of an element
\begin{align}
  c = \sum_{a \in \AA} c_a a \in A
\end{align}
is the linear combination
\begin{align}
  c^w = \sum_{a \in \AA} c_a a^w,
\end{align}
which only is an element of $A$ if $a^w \in \AA$ for all $a \in \AA$ with $c_a
\neq 0$.  Clearly, $a^w \in \AA$ if and only if $L^w \subseteq S$.  Usually
we will only consider conjugates $a^w$ with $w \in X_L^{\sharp}$.  Given $a,
a' \in \AA$ with $\iota(a) = L$ and $x \in X_L^{\sharp}$, by the definition
of the partial multiplication on $\AA$, the product $a \circ a'$ is defined
if and only if the product $a^x \circ (a')^x$ is defined, and in that case
\begin{align} \label{eq:circ^x}
  (a \circ a')^x = a^x \circ (a')^x.
\end{align}

The action of the  free monoid $S^*$ on the power set  $\PP(S)$,
together with the bijections (\ref{eq:partialbijection})
induces an
action of $S^*$ on  the set $\AA$ of all alleys as follows.   For $a = (L; s,
t, \dotsc) \in \AA$  and  $r \in S$  we set
\begin{align}
  (L; s, t, \dotsc).r := (L; s, t, \dotsc)^{\omega(L, r)}.
\end{align}
Note that $(L; s, t, \dotsc)^{\omega(L, r)} \in \AA$ since $\omega(L, r) \in X_L^{\sharp}$. 
Again, we ignore the obviously trivial edges
and define as 
 the \emph{action graph} of this action of $S^*$ on $\AA$
the graph with
vertex set $\AA$ and  edge set
$\{\edge{a}{r}{a.r} : a = (L; s, t, \dots)  \in \AA, r \in S
\setminus L\}$.

\begin{Proposition} \label{pro:a^x} Let $a, a' \in \AA$ and let $L =
  \iota(a)$.  Then $a.m = a'$ for some $m \in S^*$ if and only if $a^x = a'$
  for some $x \in X_L^{\sharp}$.
\end{Proposition}

\begin{proof} 
  Suppose first that $a' = a.m$ for some $m \in S^*$.  Let $x = \omega(L,
  m)$.  Then $\omega(L, m) \in X_L^{\sharp}$ by Proposition~\ref{pro:C-to-G}
  and $a' = a^x$.  Conversely, suppose $a' = a^x$ for some $x \in
  X_L^{\sharp}$.  Let $(L, m) \in \CC$ be a reduced expression for $(L, x)
  \in \GG$.  Then $a' = a.m$, as desired.
\end{proof}
As in Lemma~\ref{la:J.s}, if $a \in \AA$ has $\iota(a) = J$ and $L = J \cup
\{s\}$ then 
\begin{align} \label{eq:a.ss^wL}
  a.ss^{w_L} = a.
\end{align}
Hence the effects of $S^*$ on $\AA$ can be
undone and the action of the free monoid $S^*$ partitions the set of alleys
into classes of conjugate alleys.  We call such an $S^*$-orbit a \emph{street} and
denote by
\begin{align}
  [a] = [L; s, t, \dotsc] = \{a.m : m \in S^*\}
\end{align}
the class of the alley $a = (L; s, t, \dotsc) \in \AA$.
Furthermore, we let
\begin{align}
  \Psi = \{[a] : a \in \AA\}
\end{align}
denote the set of all streets.  The \emph{source} of a street $[a]$
is the shape
\begin{align}
  \iota([a]) = [\iota(a)] \in \Lambda
\end{align}
and its \emph{target} is the shape 
\begin{align}
  \tau([a]) = [\tau(a)] \in \Lambda,
\end{align}
both being well defined by Proposition~\ref{pro:a^x}.  The
\emph{length} of a street $[a]$ is 
\begin{align}
  \ell([a]) = \ell(a).
\end{align}
We write 
\begin{align}
  a \sim a'
\end{align}
if $a, a' \in \AA$ are in the same $S^*$-orbit.  Note
that in general $a^w = a'$ for some $w \in W$ does not imply $a \sim a'$.
E.g., for the element $w = s t s \in W(A_2) = \Span{s, t : s^2 = t^2 = (st)^3
  = 1}$ we have $(S;s)^w = (S^w; s^w) = (S;t)$.  But $(S;t)$ is not in the
$S^*$-orbit of $(S;s)$ since $X_S = \{\id\}$ contains no element $x$ with
$(S;s)^x = (S;t)$.

\begin{figure}[htbp]
  \centering
  \begin{tabular}{|c||c||c|c||c|}
\hline
$\emptyset\rlap{$\supset\!\!{\begin{smallmatrix}1\\2\\3\end{smallmatrix}}$}$ 
&&&& \\
\hline
\hline
\textcolor{red}{$
\begin{array}{c}
  1{;}1\rlap{$\supset_3$}\\
^2{\downarrow\uparrow}_1\\
  2{;}2\\
^3{\downarrow\uparrow}_2\\
  3{;}3\rlap{$\supset_1$}\\
\end{array}
$}
&
$
\begin{array}{c}
  1\rlap{$\supset_3$}\\
^2{\downarrow\uparrow}_1\\
  2\\
^3{\downarrow\uparrow}_2\\
  3\rlap{$\supset_1$}\\
\end{array}
$
&&& \\
\hline
\hline
\textcolor{red}{$\arraycolsep0pt
\begin{array}{ccc}
  13{;}13 & \stackrel{2}{\rightleftarrows} & 13{;}31 \\
\end{array}
$}
&
\textcolor{red}{$\arraycolsep0pt
\begin{array}{ccc}
  13{;}1 & \stackrel{2}{\rightleftarrows} & 13{;}3 \\
\end{array}
$}
&
$
\begin{array}{c}
  13\rlap{$\supset_2$}\\
\end{array}
$
&& \\
\hline
\textcolor{red}{$
\begin{array}{c}
  12{;}12\\
^3{\downarrow\uparrow}_1\\
  23{;}23\\
\end{array}
$
$
\begin{array}{c}
  12{;}21\\
^3{\downarrow\uparrow}_1\\
  23{;}32\\
\end{array}
$
}
&
$
\begin{array}{c}
  12{;}1\\
^3{\downarrow\uparrow}_1\\
  23{;}2\\
\end{array}
$
$
\begin{array}{c}
  12{;}2\\
^3{\downarrow\uparrow}_1\\
  23{;}3\\
\end{array}
$
&&
$
\begin{array}{c}
  12\\
^3{\downarrow\uparrow}_1\\
  23\\
\end{array}
$
& \\
\hline\hline
\textcolor{red}{$
\begin{array}{c}
  123{;}213\\
  123{;}123\\
  123{;}312\\
\end{array}
\begin{array}{r}
  123{;}231\\
  123{;}132\\
  123{;}321\\
\end{array}
$
}
&
$
\begin{array}{c}
\textcolor{red}{123{;}21}\\
  123{;}12\\
  123{;}31\\
\end{array}
\begin{array}{r}
\textcolor{red}{123{;}23}\\
  123{;}13\\
  123{;}32\\
\end{array}
$
&
\textcolor{red}{$
\begin{array}{c}
  123{;}2\\
\end{array}
$}
&
$
\begin{array}{c}
  123{;}3\\
\end{array}
\begin{array}{c}
  123{;}1\\
\end{array}
$
&
$
\begin{array}{c}
  123\\
\end{array}
$
\\
\hline
  \end{tabular}
  \caption{\strut Streets of $A_3$.}
  \label{tab:street-a3}
\end{figure}

\begin{Example}
  Figure~\ref{tab:street-a3} illustrates the streets for $W$ type $A_3$ as
  action graph on the set of all alleys.  Here a notation like $123;21$ is
  used as a shorthand for $(\{1,2,3\}; 2, 1) \in \AA$.  The streets are
  arranged in a square grid using their sources and targets as coordinates.
  (Some entries are highlighted in red in order to illustrate
  Proposition~\ref{pro:Delta0}.)
\end{Example}

For $L \subseteq S$ and a
subset $\BB \subseteq \AA$, the set $L \circ \BB$ is the set of all alleys $a
\in \BB$ with source $L$ and the set $\BB \circ L$ is the set of all $a \in
\BB$ with target $L$.  Moreover,
\begin{align} \label{eq:B(L)}
\BB = \coprod_{L \subseteq S} L \circ \BB 
    = \coprod_{L \subseteq S} \BB \circ L.
\end{align}
If, in particular, $\BB$ is a street
then even more can be said.

\begin{Proposition} \label{pro:alpha(L)} 
Let $a = (L; s, t, \dots) \in \AA$,
  let $\lambda = [L]$ and let $\alpha = [a]$.  Denote by
\begin{align*}
  N_a = \{x  \in N_L : a^x = a\}
\end{align*}
the stabilizer of $a$, i.e., the  stabilizer of the tuple $(s, t, \dots)$ in
$N_L$.  Then
\begin{align*}
  \alpha = \coprod_{L' \subseteq S} L' \circ \alpha,
\end{align*}
where 
\begin{enumerate}
\item $L' \circ \alpha = \emptyset$ unless $L' \in \lambda$;
\item for each $x \in X_L^{\sharp}$, the map $a \mapsto a^x$ is a bijection
  from $L \circ \alpha$ to $L^x \circ \alpha$;
\item $L \circ \alpha = a^{N_L}$ is an $N_L$-orbit of length $\Size{N_L :
    N_a}$.
\item $\Size{\alpha} = \Size{\lambda} \cdot \Size{L \circ \alpha} =
  \Size{\lambda} \cdot \Size{N_L : N_a}$.
\end{enumerate}
\end{Proposition}

\begin{proof}
  (i) and (ii) are clear.  (iii) follows from Proposition~\ref{pro:a^x}, (iv)
  from parts (i), (ii) and (iii).
\end{proof}

With the notation of Proposition~\ref{pro:alpha(L)}, we call $\Size{\lambda}$
the \emph{width} of $\alpha$ and $\Size{L \circ \alpha}$ the \emph{depth} of
$\alpha$.  Note that, by Proposition~\ref{pro:alpha(L)}(iii), the depth of a
street with source $L$ is bounded above by the order of the permutation group
induced by $N_L$ on the set $L$.

The stabilizer $N_a$ of an alley $a \in \AA$ can be described as
intersection of normalizer complements.

\begin{Theorem} \label{thm:Na-intersection}
  Let $a = (L; s_1, \dots, s_l) \in \AA$.  Then 
\begin{align*}
N_a = N_L \cap N_{L \setminus \{s_1\}}
 \cap \dots \cap N_{L \setminus \{s_1, \dots, s_l\}}.
\end{align*}
In particular, if $a, a' \in \AA$ are such that the product $a \circ a'$ is
defined, then
\begin{align*}
  N_{a \circ a'} = N_a \cap N_{a'}.
\end{align*}
\end{Theorem}

\begin{proof}
  If $\ell(a) = 0$ then $N_a = N_L$
  and there is nothing to prove.

  Otherwise, $\ell(a) = l > 0$.  Let $J = L \setminus \{s_1, \dots, s_l\}$
  and let $a' = (L; s_1, \dots, s_{l-1})$.  By induction on the length
  $\ell(a)$,
\begin{align*}
  N_{a'} = N_L \cap N_{L \setminus \{s_1\}} \cap \dots \cap
  N_{L \setminus \{s_1, \dots, s_{l-1}\}}.
\end{align*}
  Clearly, 
  \begin{align*}
    N_{a'} \cap  N_J \subseteq \{x \in N_{a'}  : J^x = J\} = \{x  \in N_{a'} :
  s_l^x  =  s_l\} =  N_a.   
  \end{align*}
Conversely, from  $N_{a'}  \subseteq  X_{LLL} \subseteq  X_{LL}
  \subseteq X_{JJ}$ it follows that 
\begin{align*}
N_a = \{x \in N_{a'} : J^x = J\} 
\subseteq  N_{a'} \cap \{x \in X_{JJ} : J^x = J\} =
  N_{a'} \cap X_{JJJ} = N_{a'} \cap N_J.
\end{align*}
Hence 
\begin{align*}
  N_{a} = N_{a'} \cap N_J = N_L
\cap N_{L \setminus \{s_1\}} \cap \dots \cap
N_{L \setminus \{s_1, \dots, s_{l-1}\}}\cap
N_{L \setminus \{s_1, \dots, s_l\}},
\end{align*}
as desired.
\end{proof}

\subsection{Products of Streets.} \label{sec:street-product}

From now on we identify a street $\alpha \in \Psi$ with the sum $\sum_{a \in
  \alpha} a$ of all its elements in~$A$.  Due to the unique factorization
property of Proposition~\ref{prop:a-category}, for each $b \in \alpha \circ
\alpha'$ there are unique factors $a \in \alpha$ and $a' \in \alpha'$ such
that $b = a \circ a'$.  Therefore, the product $\alpha \circ \alpha' \in A$
coincides with the sum over the set of products $\{a \circ a' : a \in
\alpha,\, a' \in \alpha'\}$.  In fact, the product of two streets is a sum of
streets, as the next result shows.

\begin{Theorem} \label{thm:street-product}
  Let $a = (L; s, t, \dots) \in \AA$ and $a' = (L'; s', t', \dots) \in \AA$
  be such that $L' = L \setminus \{s, t, \dotsc\}$.  
Let $D$ be  a set of double coset representatives of
$N_{a'}$ and $N_a$ in 
\begin{align*}
  N_{L'} = \coprod_{d \in D} N_{a'} d N_{a}.
\end{align*}
Then
  \begin{align*}
  [a] \circ [a'] = \sum_{d \in D} [a \circ (a')^d].
  \end{align*}
\end{Theorem}

\begin{proof}
  It suffices to show that 
  \begin{align*}
    L \circ [a] \circ [a'] = \sum_{d \in D} L \circ  [a \circ (a')^d] 
    = L \circ \sum_{d \in D} [a \circ (a')^d].  
  \end{align*}
Then, using
  Proposition~\ref{pro:alpha(L)}(i) and (ii), we can conclude that
  \begin{align*}
    J \circ [a] \circ [a'] = J
  \circ \sum_{d \in D} [a \circ (a')^d]
  \end{align*}
 for all $J \subseteq S$ and hence
  the claim follows by equation~(\ref{eq:B(L)}).

Using Proposition~\ref{pro:alpha(L)}(iii), we write
\begin{align}
  L \circ [a] \circ [a']
& = a^{N_L} \circ [a']
= \sum_{\inRCosets{n}{N_a}{N_L}} a^n \circ [a'],
\end{align} 
where $N_L/N_a$ denotes  a transversal of the  right cosets of $N_a$ in 
\begin{align}
  N_L = \coprod_{\inRCosets{n}{N_a}{N_L}} N_a n, 
\end{align}
and first calculate
$a^n \circ [a']$ for $n \in N_L$.  Then we form the sum over all these cosets
and derive the claimed formula.

Let $n \in N_L$.   From $L' \subseteq L$ follows $N_L \subseteq X_L \subseteq
  X_{L'}$, and therefore $[(a')^n] = [a']$, by Proposition~\ref{pro:a^x}.  Hence,
\begin{align}
  (L')^n \circ [a']
  &= (L')^n \circ [(a')^n]
   = (a')^{n N_{(L')^n}}
   = (a')^{n N_{L'}^n}
   = (a')^{N_{L'}n}.
\end{align}
We can thus conclude, writing $a^n = (a \circ L')^n = a^n \circ (L')^n$, that
\begin{align}
  a^n \circ [a']
  &= a^n \circ (L')^n \circ [a']
   = a^n \circ (a')^{N_{L'}n}
   = (a \circ (a')^{N_{L'}})^n,
\end{align}
which shows that
\begin{align}
  a^n \circ [a']
  &= \sum_{\inRCosets{n'}{N_{a'}}{N_{L'}}} (a \circ (a')^{n'})^n,
\end{align}
a sum over the right cosets $N_{a'} n'$ of $N_{a'}$ in $N_{L'}$. 

Let $D$ be a set of double coset representatives of 
$N_{a'}$ and $N_a$ in $N_{L'}$.
For each $d \in D$, the double coset $N_{a'} d N_a$
is a union
\begin{align}
  N_{a'} d N_a = \coprod_{\inRCosets{n''}{(N_{a'}^d \cap N_a)}{N_a}} N_{a'} dn''
\end{align}
of right cosets of $N_{a'}$, 
parametrized by the right cosets of $N_{a'}^d \cap N_a$ in $N_a$.
By Theorem~\ref{thm:Na-intersection},
\begin{align}
N_{a'}^d   \cap  N_a
= N_a \cap N_{(a')^d} 
= N_{a \circ (a')^d}.
\end{align}
Moreover, if $n' \in N_{L'}$, $d\in D$ and $n'' \in N_a$ are such that $N_{a'} n' = N_{a'} d n''$ then
\begin{align}
a \circ (a')^{n'} 
= a^{n''} \circ (a')^{d n''}
= (a \circ (a')^{d})^{n''}.
\end{align}
Therefore, if we denote $b(d) = a \circ (a')^{d}$ for $d \in D$, then
\begin{align}
\sum_{\inRCosets{n'}{N_{a'}}{N_{L'}}} (a \circ (a')^{n'})^n  
= \sum_{d \in D}
\sum_{\inRCosets{n''}{N_{b(d)}}{N_a}} b(d)^{n''n}.
\end{align}
Now $\{n''n : \inRCosets{n''}{N_{b(d)}}{N_a},\, \inRCosets{n}{N_a}{N_L}\}$ is
a transversal of the right cosets of $N_{b(d)}$ in $N_L$.  Thus a summation
over $\inRCosets{n}{N_a}{N_L}$ on both sides finally yields
\begin{align}
  L \circ [a] \circ [a']
&
=
\sum_{d \in D}
\sum_{\inRCosets{m}{N_{b(d)}}{N_L}} b(d)^m 
\\\notag&= \sum_{d \in D} b(d)^{N_L}
= \sum_{d \in D} L \circ [b(d)]
= L \circ \sum_{d \in D} [b(d)],
\end{align}
as desired.
\end{proof}

As an immediate consequence, we obtain a product formula for the depths of
streets.

\begin{Corollary}
Denote by $\dd(\alpha)$ the depth of $\alpha \in \Psi$.
  Suppose $\alpha, \alpha' \in \Psi$ and $\alpha_1, \dots, \alpha_l \in \Psi$
  are such that $\alpha \circ \alpha' = \alpha_1 + \dots + \alpha_l \neq 0$.
 Then
\begin{align*}
  \dd(\alpha) \cdot \dd(\alpha') = \dd(\alpha_1) + \dots + \dd(\alpha_l).
\end{align*}
In particular, if $\dd(\alpha) = \dd(\alpha') = 1$ then
the product $\alpha \circ \alpha'$, if defined, is a single street.
\end{Corollary}

\begin{proof}
  With the notation from (the proof of) Theorem~\ref{thm:street-product}, we
  can derive from
\begin{align}
  \Size{N_{L'}:N_{a'}} = \sum_{d \in D} \Size{N_a : N_{a'}^d \cap N_a}
\end{align}
that
\begin{align}
\dd([a]) \cdot \dd([a']) 
&= \Size{N_L:N_a} \cdot \Size{N_{L'}:N_{a'}}
= \sum_{d \in D} \Size{N_L : N_{a'}^d \cap N_a}
\\\notag&
= \sum_{d \in D} \Size{N_L : N_{a \circ (a')^d}}
= \sum_{d \in D} \dd([a \circ (a')^d]),
\end{align}
as desired.
\end{proof}

The action of $S^*$ on $\AA$ is compatible with the partial orders from
Section~\ref{sec:partial-order}, that is with taking prefixes and suffixes.
This property can be used to formulate, in
Theorem~\ref{thm:forest-product-street} below, an alternative product formula
for streets in the spirit of Proposition~\ref{pro:alley-forest-product}.

\begin{Proposition} \label{pro:action-forest}
Let $a = (L; s, t, \dots) \in \AA$ and $m \in S^*$.  Then
\begin{enumerate}
\item $\pi(a.m) = \pi(a).m$ and
\item $\sigma(a.m) =  \sigma(a).m'$, where $m' \in S^*$ is such that
$(L_s, m') \in \CC$ is a reduced expression for $(L_s, \omega(L, m)) \in \GG$.
\end{enumerate}
\end{Proposition}

\begin{proof}
  Let $d = \omega(L, m)$.  Then $a.m = a^d
  = (L^d; s^d, t^d, \dots)$.

  (i) Clearly, $\pi(a.m) = \pi(a^d) = \pi(a)^d = \pi(a).m$.

  (ii) We have $\sigma(a.m) = \sigma(a^d) = ((L_s)^d; t^d, \dots) =
  \sigma(a)^d$, since $(L^d)_{s^d} = (L_s)^d$.  From $d \in X_L \subseteq
  X_{L_s}$ and with Proposition~\ref{pro:a^x} it then follows that
  $\sigma(a.m) = \sigma(a).m'$, as desired.
\end{proof}

In other words, the preimage of any street, under $\pi$ and under $\sigma$,
is a union of streets.  The streets in the preimage of $\alpha \in \Psi$
under $\pi$ can be listed efficiently.

\begin{Proposition} \label{pro:spanning}
  Let $a  = (L; s_1, \dots, s_k) \in \AA$ and let $J = L \setminus
  \{s_1, \dots, s_k\}$.  Then
  \begin{align*}
    \pi^{-1}([a]) = \coprod_{t \in J/N_a}   [L; s_1, \dots, s_k, t],
  \end{align*}
  where $t$ ranges over a transversal of the $N_a$-orbits on $J$.
\end{Proposition}

\begin{proof}
  This follows from Proposition~\ref{pro:alpha(L)}(iii).
\end{proof}

If we define relations $\prefix$ and $\suffix$ on the set $\Psi$ of
all streets by
\begin{align}
  \alpha' \prefix \alpha \text{ if } a' \prefix a
  \text{ for some } a \in \alpha \text{ and some } a' \in \alpha'
\end{align}
and
\begin{align}
  \alpha \suffix \alpha' \text{ if } a' \suffix a
 \text{ for some } a \in \alpha \text{ and some } a' \in \alpha'
\end{align}
for $\alpha, \alpha' \in \Psi$, then 
both $\prefix$ and $\suffix$ are partial orders.
The set of streets $\Psi$ together with the two partial orders $\prefix$ and
$\suffix$ forms a rooted twofold forest with roots $\lambda \in \Lambda$.

\begin{table}[htb]
  \begin{align*}
    \begin{array}{c|ccccccccc} \hline
  n & 0 & 1 & 2 & 3 & 4 & 5 & 6 & 7 & 8 \\ \hline
A_n & 1 & 3 & 8 & 27 & 108 & 536 & 3180 & 22113 & 176175 \\
B_n &   &   & 10 & 34 & 136 & 648 & 3720 & 25186 & 196777 \\
D_n &   &   &    &    & 123 & 579 & 3417 & 23387 & 184580 \\
E_n &   &   &    &    &     &     & 3347 & 23057 & 180570 \\
F_4 &   &   &    &    & 136 \\
H_n &   &   &  8 & 30 & 120 \\
\hline
\end{array}
  \end{align*}
  \caption{The number of streets for some types of Coxeter groups.}
  \label{tab:nrStreets}
\end{table}
The streets of $W$ can be efficiently enumerated by using
Proposition~\ref{pro:spanning} to span the $\pi$-forests.  In
Table~\ref{tab:nrStreets} we list the number of streets $\Size{\Psi}$ for
some types of Coxeter groups~$W$ of small rank.

\begin{figure}[htbp]
    \centering
    \includegraphics[scale=.3]{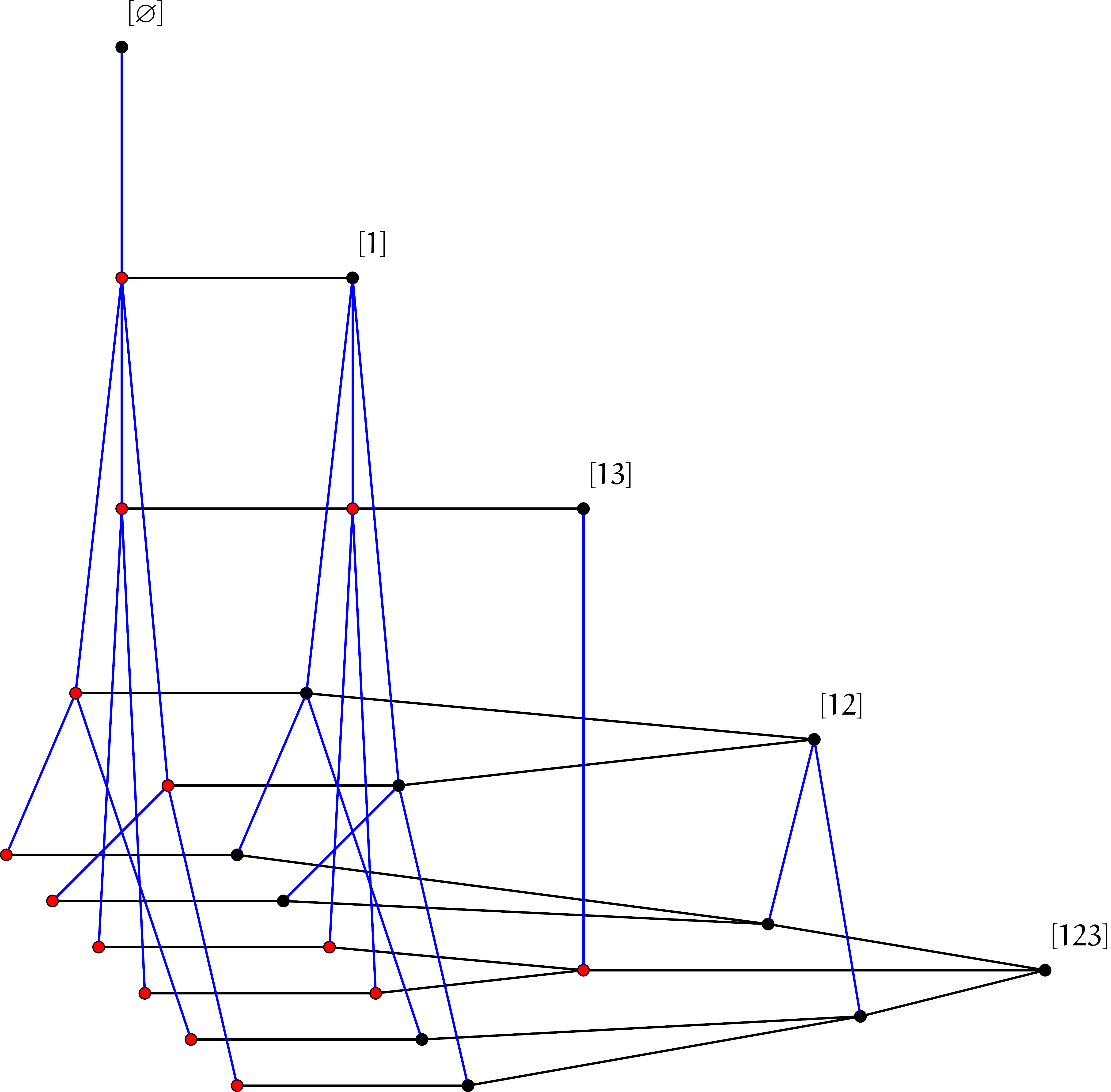}
    \caption{The rooted twofold forest for type $A_3$.}
    \label{fig:2root3}
  \end{figure}
\begin{Example}
  Figure~\ref{fig:2root3}  shows the  rooted twofold  forest for  type $A_3$.
  Here  the  $\pi$-forests  grow   horizontally  towards  the  left  and  the
  $\sigma$-forests grow  vertically downwards out  of the roots  $\lambda \in
  \Lambda$ on  the diagonal.  The other  vertices can be  identified with the
  help of Figure~\ref{tab:street-a3}.
  (Some vertices are coloured red in order to illustrate
  Proposition~\ref{pro:Delta0}.)
\end{Example}

The product of $\alpha, \alpha' \in \Psi$ can be described as the
intersection of the $\pi$-subtree spanned by $\alpha$ and the
$\sigma$-subtree spanned by $\alpha'$, provided that $\tau(\alpha) =
\iota(\alpha')$; in that case $\tau(\alpha) = \iota(\alpha')$ is the only
shape $\lambda \in \Lambda$ with $\alpha \suffix \lambda \prefix \alpha'$,
otherwise there is no such shape at all.  In any case, we have the following.

\begin{Theorem} \label{thm:forest-product-street}
Let $\alpha, \alpha' \in \Psi$.  Then
\begin{align*}
\alpha \circ \alpha' = 
 \mathop{\sum\sum}_{\substack{\lambda \in \Lambda, \alpha'' \in \Psi \\
\alpha \suffix \lambda \prefix \alpha' \\
\alpha \prefix \alpha'' \suffix \alpha'}} \alpha''.
\end{align*}
\end{Theorem}

\begin{proof}
  By definition,
  \begin{align} \label{eq:alpha-circ-alpha'}
     \alpha  \circ \alpha' = \mathop{\sum\sum}_{\substack{
a \in \alpha,\, a' \in \alpha'\\ \tau(a) = \iota(a')}}
a \circ a'.
  \end{align}
  This sum is $0$ unless $\tau(\alpha) = \iota(\alpha')$.  In that case,
  there is a unique $\lambda \in \Lambda$ with $\alpha \suffix \lambda$ and
  $\lambda \prefix \alpha'$ and
\begin{align}
  \sum_{\alpha \prefix \alpha'' \suffix \alpha'} \alpha''
&=  \sum_{\alpha \ni a \prefix a'' \suffix a' \in \alpha'} a'',
\end{align}
that is the sum over all $a'' \in \AA$ such that $a \prefix a''$ for some $a \in
\alpha$ and $a'' \suffix a'$ for some $a' \in \alpha'$, which is $\alpha
\circ \alpha'$ according to (\ref{eq:alpha-circ-alpha'}).
\end{proof}

\begin{Corollary} \label{cor:Cartan-Psi}
  Let $\alpha \in \Psi$.  Then
  \begin{enumerate}
  \item the elements $\{\alpha' \in \Psi: \alpha \prefix \alpha'\}$ form a
    basis for the right ideal $\alpha \Xi$ of $\Xi$,
  \item the elements $\{\alpha' \in \Psi: \alpha' \suffix \alpha\}$ form a
    basis for the left ideal $\Xi \alpha$ of $\Xi$.
  \end{enumerate}
  Moreover, the elements $\lambda \in \Lambda$ form a complete set of
  primitive orthogonal idempotents for the algebra $\Xi$ and the Cartan
  invariants of $\Xi$ are
  \begin{align*}
    \dim \lambda \Xi \lambda' = \# \{\alpha \in \Psi : \lambda \prefix \alpha
    \suffix \lambda'\},
  \end{align*}
  for $\lambda, \lambda' \in \Lambda$.  
\end{Corollary}

It   follows   in   particular,   that   the  matrix   $(\dim   \lambda   \Xi
\lambda')_{\lambda,\lambda'\in  \Lambda}$ of  Cartan invariants  of  $\Xi$ is
unitriangular.

Let us call $\alpha \in \Psi$ \emph{irreducible}, if $\ell(\alpha) > 0$ and
$\alpha = \alpha' \circ \alpha''$ for some $\alpha', \alpha'' \in \Psi$
implies $\ell(\alpha') \cdot \ell(\alpha'') = 0$.  Then the algebra $\Xi$ has
a presentation as a (not necessarily canonical) quiver with vertex set
$\Lambda$, edge set corresponding to the irreducible elements of $\Psi$ and
the multiplication table of $\Psi$ as its only relations.

\begin{figure}[htbp]
  \centering
  \includegraphics[scale=.3]{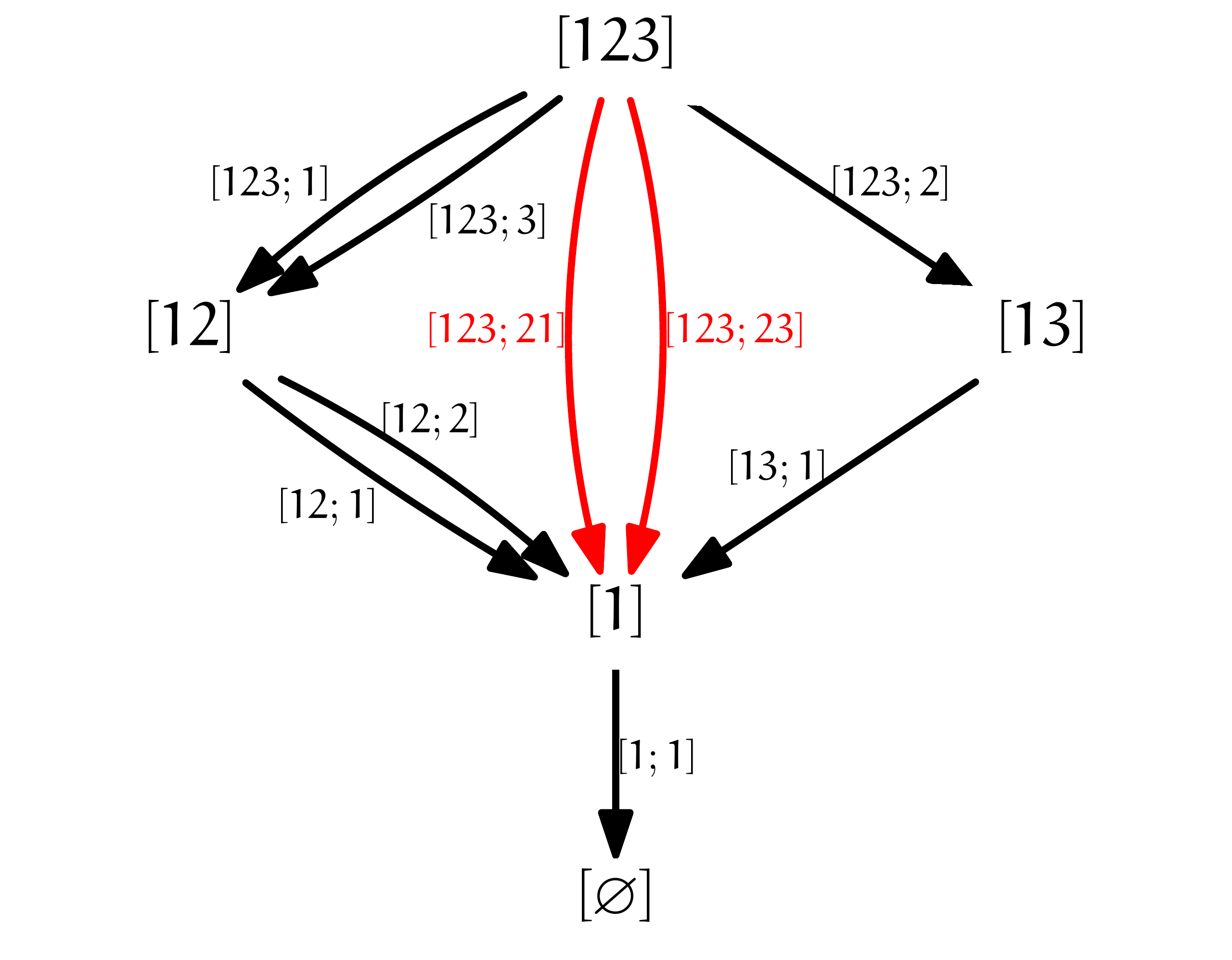}
  \caption{A quiver for $\Xi$ of type $A_3$.}
  \label{fig:street3i}
\end{figure}

\begin{Example}
  Figure~\ref{fig:street3i} illustrates the quiver for $\Xi$ in the case
  $A_3$.  There are only two irreducible streets of length greater than $1$.
  The only relevant relation in this case is $[123;2] \circ [13;1] = [123;21]
  + [123;23]$.  This relation allows us to drop one of $[123;21]$, $[123;23]$
  from the edge set.  Hence $\Xi$ of type $A_3$ is in fact a path algebra.
  The calculations for other types of Coxeter groups suggest that this is
  always the case.
\end{Example}

\begin{Conjecture}
  For any type of finite Coxeter group, the algebra $\Xi$ is a path algebra.  
\end{Conjecture}

\section{Difference Operators.} 
\label{sec:difference}

In this section a grade decreasing difference operator is introduced and
shown to eventually map $\Xi$ surjectively onto the grade $0$ component of
$A$.  In the next section, this difference operator will be used to construct
a matrix representation of $\Xi$.

For $i \geq 0$, denote by 
\begin{align}
  \EE_i =  \{\edge{a}{r}{a.r} : a = (L; s, t, \dots)  \in \AA_i, r \in S
\setminus L\}
\end{align}
the set of \emph{nontrivial  edges} of the graph of the action
of $S^*$ on  $\AA_i$.  Then the set of all edges
\begin{align}
   \EE = \coprod_{i \geq 0} \EE_i
\end{align}
is in bijection to the set of alleys $a \in \AA$ of positive length.

\begin{Proposition} \label{pro:edges}
For each $i > 0$, the map 
\begin{align*}
  (L; s, t, \dots) \mapsto \edge{(L_s; t, \dots)}{s}{(L_s; t, \dots).s}
\end{align*}
is a bijection from $\AA_i$ to $\EE_{i-1}$.
\end{Proposition}

In what follows, we will use this bijection to identify edges with elements
of $\AA$.  Note that, if $n = \Size{S}$ then $\AA_n$ has no nontrivial edges,
and therefore $\EE_n = \emptyset$.

We define the \emph{little difference} operator $\delta \colon A \to A$ as the
linear map which maps an edge $a = (L; s, t, \dotsc) \in \AA$ to the
difference of its end points, i.e.,
\begin{align}
  \label{eq:delta}
  \delta(a) =
  \begin{cases}
(L_s; t, \dotsc) - (L_s; t, \dotsc).s, & \text{if } \ell(a) > 0, \\
0, & \text{otherwise.} \\
  \end{cases}
\end{align}
Then $\delta(A_i) \subseteq A_{i-1}$ for all $i > 0$,
and hence 
\begin{align}
  \delta^i(A_i) \subseteq A_0
\end{align}
for all $i \geq 0$.
Based on this observation, we 
furthermore define the \emph{big difference} operator $\Delta \colon A \to
A_0$ as the linear map with
\begin{align}
\Delta(a) = \delta^i(a)
\end{align}
for $a \in \AA_i$.  Then, for all $a \in \AA$, we have
\begin{align}
   \Delta(a) =
   \begin{cases}
a, & \text{if } \ell(a) = 0,\\
\Delta(\delta(a)) & \text{if } \ell(a) > 0. 
   \end{cases}
 \end{align}
 We will use this latter description of $\Delta$ for inductive arguments.  It
 also follows that
\begin{align}
0 = \Delta(a) - \Delta(\delta(a)) = \Delta(a - \delta(a)),
\end{align}
whence $a - \delta(a) \in \ker \Delta$ for all $a \in A$ with $\ell(a) > 0$.

\begin{Remark}
  The  difference operator  $\delta$ does  not turn  $A$  into  a chain
  complex, or a differential graded algebra, since neither $\delta^2 = 0$ nor
  the  graded  Leibniz  Rule  $\delta(a  \circ  b)  =  \delta(a)  \circ  b  +
  (-1)^{\ell(a)} a \circ \delta(b)$ are satisfied in general. 
\end{Remark}

\begin{Proposition} \label{pro:delta.r}
  Let $a = (L; s, t, \dots) \in \AA$ and let $x \in X_L^{\sharp}$.
  Then
  \begin{enumerate}
  \item  $\delta(a^x) = \delta(a)^x$,
  \item  $\Delta(a^x) = \Delta(a)^x$.
  \end{enumerate}
\end{Proposition}

\begin{proof} (i) From $\delta(a) = \sigma(a) - \sigma(a).s$ it follows that
\begin{align*}
 \delta(a^x) &= \sigma(a^x) - \sigma(a^x).s^x
\\& = \sigma(a)^x - \sigma(a)^x.s^x
\\& = \sigma(a)^x - (\sigma(a).s)^x 
\\& = (\sigma(a) - \sigma(a).s)^x 
\\& = \delta(a)^x,
\end{align*}
as desired, since $\sigma(a^x) = \sigma(a)^x$ by
Proposition~\ref{pro:action-forest} and $(\sigma(a).s)^x = \sigma(a)^x.s^x$.

(ii) If $J \subseteq L$ then $X_L^{\sharp} \subseteq X_J^{\sharp}$.  Thus, by
(i), if $b \in \AA$ has source $\iota(b) = J$ then $\delta(b^x) =
\delta(b)^x$.  By linearity, this is even true for all $b \in \sum_{J
  \subseteq L} J \circ A$.  Noting that $a \in \sum_{J \subseteq L} J \circ
A$ implies $\delta(a) \in \sum_{J \subseteq L} J \circ A$, the result follows
by applying $\delta$ sufficiently often.
\end{proof}

The  following  graph theoretical  lemma  will help  us  to  map $\Xi$
surjectively onto $A_0 = \Q[\PP(S)]$.

\begin{Lemma} \label{la:simple-graph}
 Let $(V, E)$ be a connected simple
  graph with vertex set $V$ and edge set $E$, and let $X = \Q[V]$ be the
  vector space with basis $V$.  Then $X$ is spanned by $\sum V$ and the set
  $B = \{u - v : \{u, v\} \in E\}$.
\end{Lemma}

\begin{proof}
  Since the  graph $(V, E)$  is connected,  for any choice  of $x, y  \in V$,
  there is  a path $(v_0, v_1,  \dots, v_l)$ from $v_0  = x$ to $v_l  = y$ in
  $(V, E)$, i.e, $\{v_{i-1}, v_i\} \in E$ for $i = 1, \dots, l$, and $x - y =
  \sum_{i=1}^l (v_{i-1} -  v_i)$ is a linear combination  of elements in $B$.
  Now let $v  \in V$.  Then $\Size{V} v =  \sum V + \sum_{x \neq  v} (v - x)$
  and it follows that the space spanned by $\sum V$ and $B$ is all of $X$.
\end{proof}

Denote $\Xi_l = \Xi \cap A_l$ for $l \geq 0$.  We conclude this section
with the following important result.

\begin{Theorem}
  $A_l = \Xi_l \oplus \delta(A_{l+1})$ for all $l \geq 0$ and $A_0 =
  \bigoplus_{l \geq 0} \Delta(\Xi_l) = \Delta(\Xi)$.
\end{Theorem}

\begin{proof}
  Consider the action graph on the vertex set $\AA_l$ for some $l \geq 0$.  Its
  edge set corresponds to $\AA_{l+1}$ by Proposition~\ref{pro:edges}.  The
  connected components of this graph are the streets $\alpha \in
  \Psi_l$.  By Lemma~\ref{la:simple-graph}, the space $A_l = \Q[\AA_l]$ is
  spanned by $\Psi_l \cup \{\delta(a) : a \in \AA_{l+1}\}$.  Moreover,
  suppose an element
\begin{align}
c = \sum_{a \in \AA_l} c_a a \in A_l
\end{align}
is contained in $\Xi_l \cap \delta(A_{l+1})$.  Then $c \in \Xi_l$ implies $c_a = c_{a'}$ whenever $a \sim a'$,
whereas $c \in \delta(A_{l+1})$ implies
\begin{align}
\sum_{a \in \alpha} c_a = 0
\end{align}
for all $\alpha \in \Psi_l$.  It follows that $c_a = 0$ for all 
$a \in \AA_i$, whence $c = 0$ and
\begin{align}
A_l = \Xi_l \oplus \delta(A_{l+1}).
\end{align}
It now follows by induction from $A_{n+1} = 0$ that 
\begin{align}
A_0 = \bigoplus_{l \geq 0} \delta^l(\Xi_l) = \bigoplus_{l \geq 0} \Delta(\Xi_l) = \Delta(\Xi),
\end{align}
as desired.
\end{proof}

\section{A Matrix Representation.} \label{sec:matrix}

In this section we use the big difference operator $\Delta$ to
turn $A_0 = \Q[\PP(S)]$ into a module for the street algebra $\Xi$.
This yields an explicit matrix representation $\mu$ of~$\Xi$. 

For each $a \in A$, we define a linear map $\mu(a) \in \End A_0$
by setting
\begin{align} \label{eq:mu-alley}
  L.\mu(a) = \Delta(L \circ a)
\end{align}
for all $L \subseteq S$.

For an  alley $a  = (L; s,  t, \dots)  \in \AA$ this  means that  $L.\mu(a) =
\Delta(a)$ and  $L'.\mu(a) = 0$  for all $L'  \subseteq S$ with $L'  \neq L$.
For each  street $\alpha \in \Psi$,  we get an  endomorphism $\mu(\alpha)$ of
$A_0$ with the property that
\begin{align} \label{eq:mu-street}
  L.\mu(\alpha) 
= \sum_{a \in L \circ \alpha} \Delta(a) 
\end{align}
for all $L \subseteq S$.  The linear maps $\mu(\alpha)$, for $\alpha \in
\Psi$, have the following crucial property.

\begin{Proposition} \label{pro:delta-mu}
  $\Delta(a).\mu(\alpha') = \Delta(a \circ \alpha')$ for all $a \in A$ and
  all $\alpha' \in \Psi$.
\end{Proposition}

\begin{proof}
It suffices to consider an alley $a \in \AA$.
If $\ell(a) = 0$ then $a = \Delta(a) \subseteq S$ and
the claim is just
equation~(\ref{eq:mu-alley}).

Otherwise $\ell(a) > 0$.  Suppose that $a = (L; s, t, \dotsc)$.  Let $d =
\omega(L_s, s)$ and $L' = L \setminus \{s,t,\dotsc\}$.  Then $\delta(a) = (L_s;
t, \dotsc) - (L_s; t, \dotsc)^d$ and $\Delta(a) = \Delta(\delta(a))$.  By
induction on $\ell(a)$,
\begin{align}
\Delta((L_s; t, \dotsc)).\mu(\alpha') 
&= \Delta((L_s; t, \dotsc) \circ \alpha') \\\notag
&= \Delta\Bigl(\sum_{a' \in L' \circ \alpha'} (L_s; t, \dotsc) \circ a'\Bigr) \\\notag
&= \sum_{a' \in L' \circ \alpha'} \Delta\bigl( (L_s; t, \dotsc) \circ a'\bigr).
\end{align}
Similarly,
\begin{align}
\Delta((L_s; t, \dotsc)^d).\mu(\alpha') 
&=\Delta((L_s; t, \dotsc)^d \circ \alpha') \\\notag
&= \Delta\Bigl(\sum_{a' \in L' \circ \alpha'} (L_s; t, \dotsc)^d \circ (a')^d\Bigr) \\\notag
&= \sum_{a' \in L' \circ \alpha'}  \Delta\bigl( ((L_s; t, \dotsc) \circ a')^d\bigr).
\end{align}
Hence
\begin{align}
  \Delta(a).\mu(\alpha') 
&= \left(\Delta(\mathstrut(L_s; t, \dotsc)) - \Delta((L_s; t, \dotsc)^d)\right).\mu(\alpha') \\\notag
&= \sum_{a' \in L' \circ \alpha'} \Delta\bigl( (L_s; t, \dotsc) \circ a'\bigr) - \Delta\bigl( ((L_s; t, \dotsc) \circ a')^d\bigr) \\\notag
&= \sum_{a' \in L' \circ \alpha'} \Delta(\delta(a \circ a')) \\\notag
&=  \Delta(a \circ \alpha'),
\end{align}
using the facts that $(L_s \setminus \{t, \dots\})^d = L_s^d \setminus \{t^d, \dots\}$ and that $(a'')^d \circ (a')^d = (a'' \circ a')^d$ by equation~(\ref{eq:circ^x}).
\end{proof}

We denote the restriction of the linear map $\mu$ to the subalgebra $\Xi$ of
$A$ again by~$\mu$.

\begin{Theorem}
  The map $\mu \colon \Xi \to \End A_0$ defined by 
  \begin{align*}
    L.\mu(\alpha) = \Delta(L  \circ \alpha)
  \end{align*}
  is a homomorphism of algebras.
\end{Theorem}

\begin{proof}
Let $a \in A$ and $\alpha' \in \Psi$.
By Proposition~\ref{pro:delta-mu},   we have
\begin{align}
  L.\mu(a) \mu(\alpha') 
&= \Delta(L \circ a). \mu(\alpha')
= \Delta(L \circ a \circ \alpha') 
= L.\mu(a \circ \alpha'),
  \end{align}
  for all $L \subseteq S$.  It follows that $\mu(a \circ \alpha') = \mu(a)
  \mu(\alpha')$ for all $a \in A$ and all $\alpha' \in \Psi$ and thus, in
  particular, $\mu(\alpha \circ \alpha') = \mu(\alpha) \mu(\alpha')$ for all
  $\alpha, \alpha' \in \Psi$.
\end{proof}

\section{More about Descents.}
\label{sec:more-descents}

In this section we identify the $\Xi$-module $A_0$ with the descent algebra
$\Sigma(W)$ and show that the linear maps $\mu(\alpha)$ for $\alpha \in \Psi$
are endomorphisms of $\Sigma(W)$.  It follows that $\Sigma(W)$ is
anti-isomorphic to $\Xi/\ker \mu = \Xi/\ker \Delta$.  This gives us the
desired presentation of $\Sigma(W)$ as a quiver with relations.

We first take a closer look at certain of the sets $X_{JKL}$ from
equation~(\ref{eq:xJKL}).

\begin{Proposition} \label{pro:xmjj-xmll} Let $s \in S$ and $J, K, L, M
  \subseteq S$ be such that $J \cup \{s\} = L$ and $K = J^d$, where $d =
  \omega(J, s)$.
Then
\begin{enumerate}
\item $X_{MJJ} \cap X_{MKK} = X_{MLL}$;
\item $(X_{MJJ} \setminus X_{MLL}) d = X_{MKK} \setminus X_{MLL}$.
\end{enumerate}
\end{Proposition}

\begin{proof}
  (i) First note that $J \subseteq L$ implies $X_L \subseteq X_J$ and
  \begin{align}
  X_{MLL} = \{x \in X_M \cap X_L^{-1} : M^x \supseteq L\} \subseteq \{x \in
  X_M \cap X_J^{-1} : M^x \supseteq J\} = X_{MJJ}.
  \end{align}
  Similarly, $X_{MLL} \subseteq X_{MKK}$.  Conversely, let $x \in
  X_{MJJ} \cap X_{MKK}$.  Then $x \in X_M$, and $x \in X_J^{-1} \cap X_K^{-1} =
  X_L^{-1}$, and $M^x \supseteq J \cup K = L$ whence
  $x \in X_{MLL}$.

  (ii) Let $x \in X_{MJJ} \setminus  X_{MLL}$. By (i) and by symmetry it
  suffices to show that $xd  \subseteq X_{MKK}$.  We have $xd \in
  X_J^{-1}  d  =  X_K^{-1}$  by Proposition~\ref{pro:xJ}(ii) and $M^{x  d}  \cap  K  =
  (M^x \cap  J)^d = J^d  = K$.  It  remains to show that  $xd \in
  X_M$.  We distinguish two cases.

  If $x \in X_L^{-1}$ then $x \in X_{ML}$ and $x \in X_{MJJ} \setminus
  X_{MLL}$ implies $M^x \cap L = J$.  Hence, using the Mackey decomposition
  of Proposition \ref{pro:xJ}(iii),
  \begin{align}
  X_M = \coprod_{b \in X_{ML}} b X_{M^b \cap L}^L \supseteq x X_J^L
  \end{align}
  and $xd \in X_M$, since $d = w_J w_L \in X_J^L$.

  Otherwise, $x \notin X_L^{-1}$.  Then $s$ is a prefix of $x^{-1}$ and since
  $x^{-1} \in X_J$ and $J^{x^{-1}} \subseteq M$, Lemma~\ref{la:sharp}(iii)
  implies that $d$ is a prefix of $x^{-1}$.  Therefore, $xd \in X_M$ as $xd$
  then is a prefix of $x \in X_M$.
\end{proof}

Following Bergeron, Bergeron, Howlett and Taylor \cite{BeBeHoTa92}, we
define numbers $m_{KL}$, for $K, L \subseteq S$, as
\begin{align}
  m_{KL} &= \sum_{J \sim L} a_{JKL} = 
\begin{cases}
\Size{X_K \cap X_L^{\smash{\sharp}}}, 
& \text{if } L \subseteq K, \\
 0,  & \text{otherwise.}
\end{cases}
\end{align}
Then,  for   a  suitable  ordering  of   the  subsets  of   $S$,  the  matrix
$(m_{KL})_{K,L  \subseteq  S}$  is  lower triangular  with  nonzero  diagonal
entries
\begin{align}
m_{KK} = \sum_{J \sim K} a_{JKK} = \#[K] \cdot \Size{N_K},
\end{align}
for $K  \subseteq S$,  and thus has  an inverse $(n_{JK})_{J,K  \subseteq S}$
over $\Q$.

\begin{Proposition}  \label{pro:eJ} 
  The  elements $e_J \in \Sigma(W)$, defined for $J  \subseteq S$ as
  \begin{align} \label{eq:e-as-x}
  e_J = \sum_K n_{JK}\, x_K,
  \end{align}
  form a basis of $\Sigma(W)$ with
  \begin{align} \label{eq:e.x}
    e_J\, x_M &= \sum_{K \sim J} a_{JMK}\, e_K
  \end{align}
 for all $J, M \subseteq S$.
  Moreover, the elements $e_{\lambda}$, defined for $\lambda \in \Lambda$ as
\begin{align}
e_{\lambda} = \sum_{L \in \lambda} e_L,
\end{align}
form a complete set of primitive orthogonal idempotents of $\Sigma(W)$ with
$e_{\lambda} e_M = e_M$ if $M \in \lambda$, and $e_{\lambda} e_M = 0$,
otherwise.
\end{Proposition}

\begin{proof}\relax
 \cite[Theorem 7.8 and Proposition 7.11]{BeBeHoTa92}.
\end{proof}

We further define, for each alley $a = (L; s, t, \dotsc) \in \AA$, an element
$f_a \in \Sigma(W)$ as
\begin{align}
  f_a &= 
  \begin{cases}
e_L, & \text{ if } \ell(a) = 0, \\
  f_{\sigma(a)} - f_{\sigma(a).s}, & \text{ if } \ell(a) > 0.
  \end{cases}
\end{align}
More generally, for a linear combination $c = \sum_{a \in \AA} c_a a \in A$,
we define
\begin{align}
  f_c = \sum_{a \in \AA} c_a f_a. 
\end{align}
Then $f_a = f_{\delta(a)}$ for all  $a \in \AA$ and $f_c = f_{\Delta(c)}$ for
all $c \in A$.  In particular, for  each street $\alpha \in \Psi$, we have an
element
\begin{align}
  f_{\alpha} &= \sum_{a \in \alpha} f_a,
\end{align}
and for each $L \subseteq S$ we have
\begin{align}
  f_{L \circ \alpha} &= \sum_{a \in L \circ \alpha} f_a,
\end{align}
If we identify $A_0$ with $\Sigma(W)$ by setting $(L; \emptyset) = e_L$,
we even have $f_c = \Delta(c)$ for all $c \in A$.

By Proposition~\ref{pro:eJ}, right multiplication by $x_M$ maps $e_J$ to a
sum of conjugates of $e_J$.  The following key result
generalizes this property to the elements $f_a$.  Recall from
Section~\ref{sec:streets} that, if $a = (L; s, t, \dots) \in \AA$ and $x \in
X_L^{\sharp}$, then $a^x = (L^x; s^x, t^x, \dots) \in \AA$.

\begin{Theorem} \label{thm:faxm}
  Let $L, M \subseteq S$ and $a = (L; s, t, \dots) \in \AA$.  Then
\begin{align*}
f_a x_M = \sum_{x \in X_{MLL}^{-1}} f_{a^x}.
\end{align*}
In particular, $f_a x_M = 0$ unless $L$ is contained in a conjugate of $M$.
\end{Theorem}

\begin{proof}
  If $\ell(a) = 0$ then $f_a = e_L$ and by equation~(\ref{eq:e.x}),
  \begin{align}
  e_L x_M
  = \sum_{K \sim L} a_{LMK} e_K
  = \sum_{x \in X_{LM} : L^x \subseteq M} e_{L^x \cap M}
  = \sum_{x \in X_{MLL}^{-1}}
  e_{L^x},
  \end{align}
  since $\{x \in X_{LM} : L^x \subseteq M\} = 
  \{x \in X_{ML}^{-1} : M^{x^{-1}} \cap L = L\} = X_{MLL}^{-1}$.

  Otherwise $\ell(a) > 0$ and $f_a = f_{(J; t, \dots)} -
  f_{(K; t', \dots)}$, where $J  = L_s$, $d = \omega(J, s)$, $K = J^d$ and
  $t' = t^d$.  By induction on $\ell(a)$ and by
  Proposition~\ref{pro:xmjj-xmll}, we have
  \begin{align}
  f_{(J; t, \dots)} x_M &= \sum_{x \in X_{MJJ}^{-1}} f_{(J; t, \dots)^x}
\\\notag&
= \sum_{x \in X_{MLL}^{-1}} f_{(J; t, \dots)^x}
+ \sum_{x \in X_{MJJ}^{-1} \setminus X_{MLL}^{-1}} f_{(J; t, \dots)^x}
  \end{align}
  and
  \begin{align}
  f_{(K; t', \dots)} x_M &= \sum_{y \in X_{MKK}^{-1}} f_{(K; t', \dots)^y}
\\\notag&
= \sum_{y \in X_{MLL}^{-1}} f_{(K; t', \dots)^y}
+ \sum_{x \in X_{MJJ}^{-1} \setminus X_{MLL}^{-1}} f_{(J; t, \dots)^x}.
  \end{align}
It follows with Proposition~\ref{pro:delta.r} that
  \begin{align}
    f_a x_M &= (f_{(J; t, \dots)} -  f_{(K; t', \dots)}) x_M 
\\\notag&
= \sum_{x \in X_{MLL}^{-1}} (f_{(J; t, \dots)^x} -  f_{(K; t', \dots)^x}) 
= \sum_{x \in X_{MLL}^{-1}} f_{a^x},
  \end{align}
as desired.
\end{proof}

\begin{Corollary} \label{cor:faxm} Let $a = (L; s, t, \dots) \in \AA$, let
  $\alpha = [a]$, let $\lambda = [L]$ and let $M \subseteq S$.   If $M
  \in \lambda$ then
\begin{align*}
 \text{\textup{(i)} } f_a e_M &=  \frac1{\Size{\alpha}} f_{M \circ \alpha}, &
 \text{\textup{(ii)} }  f_{L \circ \alpha} e_M &= \frac1{\Size{\lambda}} f_{M \circ \alpha}, &
 \text{\textup{(iii)} }  f_{\alpha} e_M &= f_{M \circ \alpha}.
\end{align*}
Otherwise, $f_{\alpha} e_M = 0$.
\end{Corollary}

\begin{proof}
  Suppose first that $M \in \lambda$.  (i) From equation~(\ref{eq:e-as-x}),
  we have
\begin{align}
f_a e_M = f_a \sum_J n_{MJ} x_J = \sum_{J \subseteq M} n_{MJ} f_a x_J =
n_{MM} f_a x_M = \frac1{m_{MM}} f_a x_M
\end{align}
since, by Theorem~\ref{thm:faxm}, $f_a x_J = 0$ unless $L$ is contained in a
conjugate of $J \subseteq M$.  Moreover, $m_{MM} = \Size{N_L} \Size{\lambda}$
and by Theorem~\ref{thm:faxm}, $f_a x_M = \Size{N_a} f_{M \circ \alpha}$
since $M \circ \alpha$ is the $N_M$-orbit of $a \in \alpha$ and $N_a$ is its
stabilizer.   The claim now follows from $\Size{\alpha} =
\Size{N_L: N_a} \Size{\lambda}$, see Proposition~\ref{pro:alpha(L)}(iv).

(ii) and (iii) follow easily from (i).

Now suppose that $M \notin \lambda$.  With $f_{\lambda} = \sum_{L' \in \lambda}
e_{L'}$, we have
\begin{align}
   f_{\alpha} f_{\lambda} 
= \sum_{L' \in \lambda} f_{\alpha} e_{L'}
= \sum_{L' \in \lambda} f_{L' \circ \alpha} 
= f_{\alpha}.
\end{align}
It follows that, for $M \notin \lambda$,
\begin{align}
f_{\alpha} e_M = f_{\alpha} f_{\lambda} e_M = 0,
\end{align}
since $f_{\lambda}  e_M  =  0$ unless  $M  \in  \lambda$.
\end{proof}

We identify $A_0$ with $\Sigma(W)$ by setting 
\begin{align}
  L = (L; \emptyset) = e_L
\end{align}
for all $L \subseteq S$ and can now formulate the main result of this article.

\begin{Theorem} \label{thm:main} The linear map $\Delta \colon \Xi \to
  \Sigma(W)$ defined by $\alpha \mapsto \Delta(\alpha) = f_{\alpha}$ for
  $\alpha \in \Psi$ is a surjective anti-homomorphism of algebras which
  induces a bijection between the complete set of primitive orthogonal
  idempotents $\lambda \in \Lambda$ of $\Xi$ and the complete set of
  primitive orthogonal idempotents $e_{\lambda}$, $\lambda \in \Lambda$, of
  $\Sigma(W)$.  
\end{Theorem}

\begin{proof}
  Let $\alpha, \alpha' \in \Psi$ and let $M \subseteq S$.  
Then 
\begin{align}
f_{\alpha} e_M &= f_{M
    \circ \alpha} = \Delta(M \circ \alpha) = M.\mu(\alpha) =
  e_M.\mu(\alpha).
\end{align}
It follows that
  \begin{align}
    f_{\alpha'} f_{\alpha} e_M
&= f_{\alpha'} (e_M.\mu(\alpha))
= (e_M.\mu(\alpha)).\mu(\alpha')
\\\notag&= e_M.(\mu(\alpha)\mu(\alpha'))
= e_M.\mu(\alpha \circ \alpha')
= f_{\alpha \circ \alpha'} e_M.
  \end{align}
Hence 
\begin{align}
  f_{\alpha'} f_{\alpha} = f_{\alpha \circ \alpha'},
\end{align}
as desired.
\end{proof}

Note that the linear map defined by $a \mapsto f_a$ in general is not an
algebra homomorphism from $A$ to $\Sigma(W)$: the product of $(J, \emptyset)$
and $(K, \emptyset)$ in $\Q\AA$ is zero unless $J = K$ while $e_J e_K \neq 0$
if $J \sim K$ by Proposition~\ref{pro:eJ}.

As an immediate consequence of Theorem~\ref{thm:main}, we derive some
properties of the Cartan matrix of $\Sigma(W)$ from the Cartan matrix
of~$\Xi$.

\begin{Corollary} \label{cor:cartan} The Cartan invariants of $\Sigma(W)$ are
  given by the dimensions of the subspaces
  \begin{align*}
  e_{\lambda'} \Sigma(W) e_{\lambda} 
    =  \Span{f_{\alpha} : \lambda \prefix \alpha \suffix \lambda'}_{\Q}
  \end{align*}
  of $\Sigma(W)$, for all $\lambda, \lambda' \in \Lambda$.  Hence the matrix
  $(\dim e_{\lambda'} \Sigma(W) e_{\lambda})_{\lambda', \lambda \in \Lambda}$
  of Cartan invariants of $\Sigma(W)$ is unitriangular.
\end{Corollary}

\begin{proof}
  By Corollary~\ref{cor:Cartan-Psi}, the set $\{\alpha \in \Psi : \lambda
  \prefix \alpha \suffix \lambda'\}$ forms a basis of the subspace $\lambda
  \Xi \lambda'$ of $\Xi$.  Under the anti-homomorphism $\Delta$, this set is
  mapped to $\{f_{\alpha} : \lambda \prefix \alpha \suffix \lambda'\}$ which
  therefore spans the subspace
\begin{align*}
  e_{\lambda'} \Sigma(W) e_{\lambda} = \Delta(\lambda \Xi \lambda')
\end{align*}
of $\Sigma(W)$.
\end{proof}

\section{The Quiver of the Descent Algebra.}
\label{sec:properties}

Denote by $\mathbf{Q} = (\mathbf{V}, \mathbf{E})$ the quiver of $\Sigma(W)$.
This is a graph with vertex set $\Vb$ corresponding to the shapes $\Lambda$
of $W$, which, by Proposition~\ref{pro:eJ}, label the complete set of
primitive orthogonal idempotents $f_{\lambda} = e_{\lambda}$, $\lambda \in
\Lambda$, of $\Sigma(W)$, and edge set $\Eb$ consisting of $\dim
e_{\lambda'}(\Rad \Sigma(W)/\Rad^2 \Sigma(W)) e_{\lambda}$ edges from
$\lambda'$ to $\lambda$ for all $\lambda, \lambda' \in \Lambda$.

We denote by $\leq$ the partial order induced on $\Lambda$ by subset
inclusion and by $\lessdot$ the cover relation of this partial order,
i.e.,  given $\lambda, \lambda' \in \Lambda$ we write
\begin{align}
  \lambda' \leq \lambda 
\end{align}
if $L' \subseteq L$ for some 
$L \in \lambda$, $L' \in \lambda'$, and
\begin{align}
  \lambda' \lessdot \lambda
\end{align}
if $\lambda' < \lambda$ and there is no $\rho \in \Lambda$ with $\lambda' <
\rho < \lambda$.

It follows from Corollary~\ref{cor:cartan} that, if there is an edge $\eb$
from $\lambda'$ to $\lambda$ in the quiver $\Qb$ then $\lambda' < \lambda$.
Some further properties of $\Sigma(W)$ and its quiver $\Qb$ follow easily
from the description of $\Sigma(W)$ as anti-homomorphic image of the streets
algebra $\Xi$.  We can, for example, find some streets in $\ker \Delta$.
Recall from Proposition~\ref{pro:edges} that an alley $a = (L; s,t, \dots)$
corresponds to an edge from $\sigma(a)$ to $\sigma(a).s$ in the action graph
on~$\AA$.  If this edge is a loop then $a$ and its $S^*$-orbit $[a]$ lie in
$\ker \Delta$.  In fact, each street $\alpha$ which ends in $[a]$ then lies in
$\ker \Delta$, as the next result shows.

\begin{Proposition} 
  \label{pro:Delta0} 
  Let  $a =  (L;  s,t,  \dots) \in  \AA$  be such  that  $\ell(a)  > 0$.   If
  $\sigma(a) = \sigma(a).s$  then $f_{\alpha} = 0$ for  all $\alpha \in \Psi$
  such that $\alpha \suffix [a]$.
\end{Proposition}

\begin{proof}
  We first consider the case $\alpha = [a]$.  From $\sigma(a) = \sigma(a).s$
  follows $\delta(a) = 0$ and thus $\Delta(a) = \Delta(\delta(a)) = \Delta(0)
  = 0$.  Moreover, by Proposition~\ref{pro:delta.r}, we have $\Delta(a.m) =
  0$ for all $m \in S^*$.  Hence $f_{\alpha} = \Delta(\alpha) = \sum_{a \in
    \alpha} \Delta(a) = 0$.

  Now suppose $\alpha = [a'] \in \Psi$ is such that $\Delta(\sigma(a')) = 0$.
  Then $\Delta(a') = 
\Delta(\delta(a')) = 
0$ and it  follows as before that $f_{\alpha} = 0$.  The
  claim for all $\alpha \suffix [a]$ then follows by induction.
\end{proof}

\begin{Example}
  In Figure~\ref{tab:street-a3} and  Figure~\ref{fig:2root3}, for $W$ of type
  $A_3$, all those $\alpha \in \Psi$ with $f_{\alpha} = 0$ are highlighted in
  red.  In  this example, all cases of  $f_{\alpha} = 0$ can  be explained by
  Proposition~\ref{pro:Delta0}
\end{Example}

A street $\alpha$ may be in $\ker \Delta$ for other reasons.
By equation~(\ref{eq:a.ss^wL}), an alley $a = (L; s, t, u, \dotsc) \in \AA$
of length $l(a) > 0$ corresponding to the edge $\edge{(L_s; t, u,
  \dots)}{s}{(L_s; t, u, \dots).s}$ of the action graph has a reverse edge
\begin{align}
  \edge{(L_s; t, u, \dots).s}{s'}{(L_s; t, u, \dots)}
\end{align}
going in the opposite direction, where  $s' = s^{w_L}$.  To this reverse edge
corresponds the \emph{reverse alley} $\overline{a}$ of $a$ which we accordingly
define as
\begin{align}
  \overline{a} = (L; s^{w_L}, t^d, u^d, \dotsc),
\end{align}
where $d = \omega(L_s, s)$.  Clearly,
\begin{align}
  \delta(\overline{a}) = - \delta(a)
\end{align}
and
\begin{align}
  {\overline{a}}^x = \overline{a^x}
\end{align}
for all $x \in X_L^{\sharp}$.  From this, the following properties are
obvious.

\begin{Lemma} \label{la:reverse0} Let $\alpha \in \Psi$.  Then
  \begin{enumerate}
  \item $\overline{\alpha} = \{\overline{a} : a \in \alpha\} \in \Psi$ with
    $\iota(\overline{\alpha}) = \iota(\alpha)$ and $\tau(\overline{\alpha}) =
    \tau(\alpha)$;
  \item $f_{\overline{\alpha}} = - f_{\alpha}$;
  \item $f_{\alpha} = 0$ if  $\overline{\alpha} = \alpha$.
  \end{enumerate}
\end{Lemma}

Lemma~\ref{la:reverse0}(iii) has consequences for the images of streets
$\alpha$ of length $\ell(\alpha) = 1$ in $\Sigma(W)$.

\begin{Proposition} 
  \label{pro:fLs=0} 
  Let $s \in L \subseteq S$.  If there exists a subset $M$ with $L \subseteq
  M \subseteq S$ such that $(L; s)^{w_M} = (L; s)$ then $f_{[L; s]} =
  0$.  In particular, if the longest element $w_0$ is central in $W$ then
  $f_{[L; s]} = 0$ for all $L \subseteq S$ and $s \in L$.
\end{Proposition}

\begin{proof}
  The reverse of the alley $(L; s)$  is $(L; s^{w_L})$.  Clearly $x = w_L w_M
  \in  X_L^{\sharp}$.  The  claim now  follows  with Lemma~\ref{la:reverse0},
  since $(L; s^{w_L})^x = (L^{w_M}; s^{w_M}) = (L; s)$.
\end{proof}

\subsection{Reduction to Irreducible Finite Coxeter Groups.}
\label{sec:reduction}

It has been shown~\cite[Proposition 3.2]{BeBeHoTa92}, that if there are
subsets $J, K \subseteq S$ such that $W$ is the direct product $W_J \times
W_K$ then the descent algebra $\Sigma(W)$ is the tensor product of
$\Sigma(W_J)$ and $\Sigma(W_K)$.  The quiver $\Qb = (\Vb, \Eb)$ of
$\Sigma(W)$ is then the direct product of the quiver $\Qb_J = (\Vb_J, \Eb_J)$
of $\Sigma(W_J)$ and the quiver $\Qb_K = (\Vb_K, \Eb_K)$ of $\Sigma(W_K)$.
It has vertex set $\Vb = \Vb_J \times \Vb_K$ and edges $(x', y') \to (x, y)$
if $x' = x \in \Vb_J$ and $y' \to y$ is an edge in $\Eb_K$ or if $y' = y \in
\Vb_K$ and $x' \to x$ is an edge in $\Eb_J$.  Moreover, in every square
arising from the product of an edge $x' \to x$ in $\Eb_J$ and an edge $y' \to
y$ in $\Eb_K$, the relation
\begin{align}
  \left((x', y') \to (x', y) \to (x, y)\right) 
= \left((x', y') \to (x, y') \to (x, y)\right)
\end{align}
holds in $\Sigma(W)$.

The problem of finding a quiver presentation for $\Sigma(W)$ is thus reduced
to irreducible Coxeter groups and we will, for the remainder of this section,
assume that $W$ is an irreducible finite Coxeter group.  The following
property of maximal parabolic subgroups is then easily verified in a
case-by-case analysis.

\begin{Proposition} 
  \label{pro:notH3} 
  Suppose that $W$ is an irreducible
  finite Coxeter group.  Let $L \subseteq S$ and let $s, s' \in L$ be such
  that $L_s \sim L_{s'}$ in $W$.  Then there exists an $x \in N_W(W_L)$ such
  that $L_s^x = L_{s'}$, unless $W$ is of type $H_3$ and $W_L$ of type $A_1
  \times A_1$.
\end{Proposition}

\begin{Corollary} 
  \label{cor:notH3}
  Suppose that $W$ is an irreducible finite Coxeter group not of type $H_3$.
  Let $\lambda, \lambda' \in \Lambda$ be such that $\lambda' \lessdot
  \lambda$ and let $\alpha, \alpha' \in \Psi$
have common source $\lambda$ and target $\lambda'$.
  Then $\alpha' = \alpha$ or $\alpha' = \overline{\alpha}$.
\end{Corollary}

\begin{proof}
  Let $s,s' \in L \subseteq S$ be such that $L \in \lambda$ and $L_s, L_{s'}
  \in \lambda'$.  Then, by Proposition~\ref{pro:notH3}, $L_s^x = L_{s'}$ for
  some $x \in N_W(W_L) = W_L \rtimes N_L$.  Let $J = L_s$ and write $x = y
  \cdot z$ for $y \in X_J^L$ and $z \in N_L$.  Then $J$ is a maximal subset
  of $L$ and thus either $y = \id$ or $y = w_J w_L$, by Lemma~\ref{la:sharp}.
  With Proposition~\ref{pro:alpha(L)}(iii) it now follows that either $(L;s')
  \sim (L; s)$ or $(L;s') \sim (L;s^{w_L})$.
\end{proof}

\begin{Remark}
  A direct inspection shows that, if $\lambda,\lambda' \in \Lambda$ are such
  that $\lambda' \lessdot \lambda$ then $\# \{\alpha \in \Psi : \lambda
  \prefix \alpha \suffix \lambda'\} \in \{1, 2\}$ holds for $W$ of type $H_3$
  as well.
\end{Remark}

Each cover relation gives rise to at most one edge in $\Qb$.

\begin{Theorem} 
  \label{thm:edge1} 
  Let $\lambda, \lambda' \in \Lambda$ be such that $\lambda' \lessdot
  \lambda$ and denote by $n_{\lambda\lambda'}$ the number of edges from
  $\lambda'$ to $\lambda$ in $\Qb$.  If there are $s \in L \subseteq M
  \subseteq S$ such that $L \in \lambda$, $L_s \in \lambda'$ and $(L;
  s)^{w_M} = (L;s)$ then $n_{\lambda\lambda'} = 0$.  Otherwise
  $n_{\lambda\lambda'} \leq 1$.
\end{Theorem}

\begin{proof}
  The number of edges from $\lambda'$ to $\lambda$ is given by 
\begin{align}
n_{\lambda\lambda'} 
= \dim  e_{\lambda'} (\Rad \Sigma(W)/\Rad^2 \Sigma(W)) e_{\lambda}
= \dim  e_{\lambda'} \Sigma(W) e_{\lambda}, 
\end{align}
since clearly $e_{\lambda'} \Sigma(W) e_{\lambda} \leq \Rad \Sigma(W)$ and
$e_{\lambda'} \Sigma(W) e_{\lambda} \cap \Rad^2 \Sigma(W) = 0$.  And by
Corollary~\ref{cor:cartan}, $e_{\lambda'} \Sigma(W) e_{\lambda}$ is spanned
by $\{f_{\alpha} : \lambda \prefix \alpha \suffix \lambda'\}$.  

If $W$ is of type $H_3$ then the longest element $w_0$ is central in $W$ and
it follows from Proposition~\ref{pro:fLs=0} that $n_{\lambda\lambda'} = 0$.

Otherwise, by Corollary~\ref{cor:notH3} and Lemma~\ref{la:reverse0}(ii), the
subspace $e_{\lambda'} \Sigma(W) e_{\lambda}$ is spanned by a single element
$f_{\alpha} = -f_{\overline{\alpha}}$, which, by Proposition~\ref{pro:fLs=0},
is $0$ if there are $s \in L \subseteq M \subseteq S$ such that $L \in
\lambda$, $L_s \in \lambda'$ and $(L; s)^{w_M} = (L;s)$.
\end{proof}

  \begin{figure}[htbp]
    \centering
    \includegraphics[scale=.3]{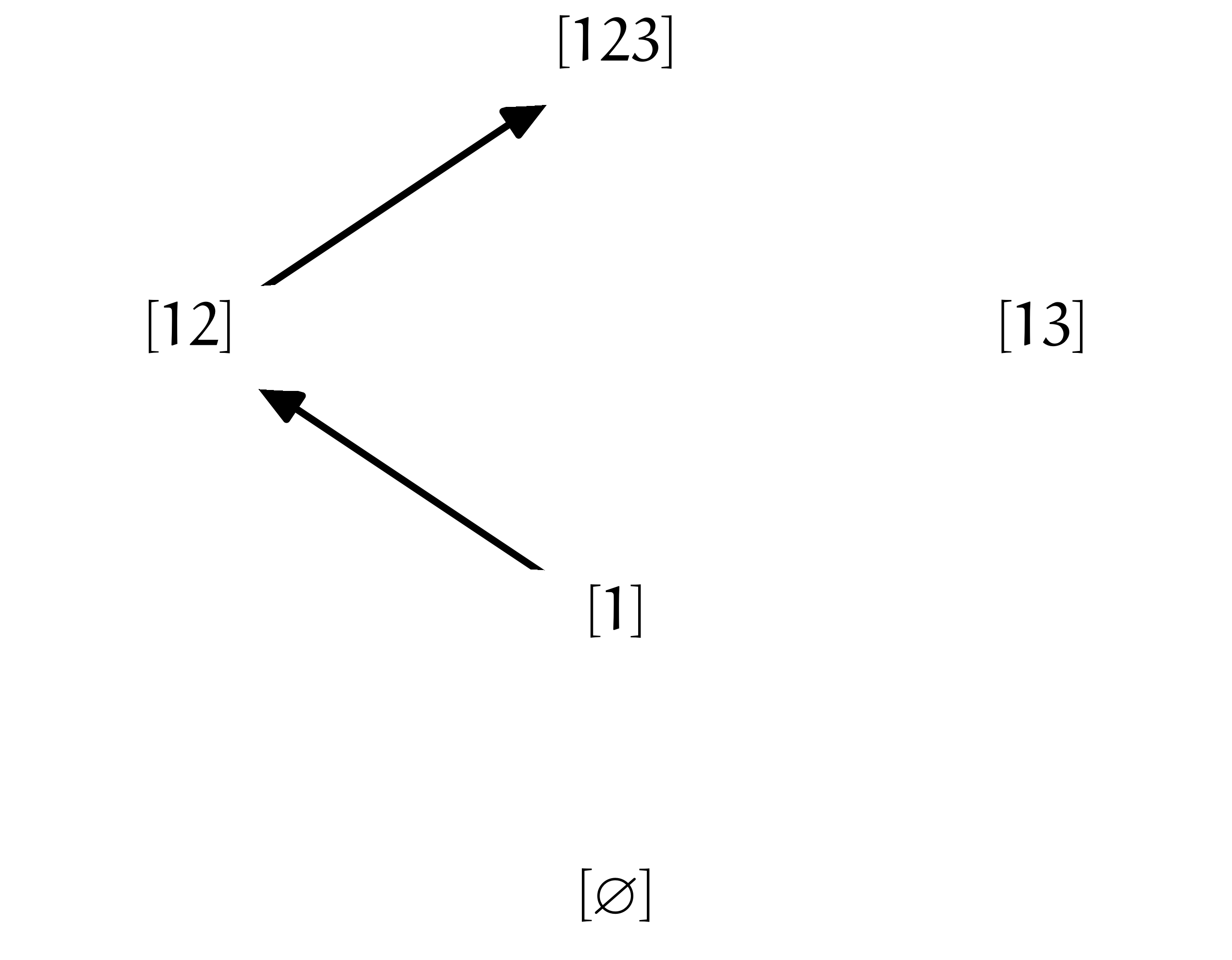}
    \caption{The quiver of $\Sigma(W)$ of type $A_3$.  There are no relations.}
    \label{fig:descent3}
  \end{figure}
\begin{Example}
  Figure~\ref{fig:descent3} illustrates the case $A_3$.  The vertices of the
  quiver are the shapes $\Lambda$ of $W$, which in this case correspond to
  the partitions of~$4$.  There are only two edges, one mapping to
  $f_{[12;1]}$ and the other to $f_{[123;1]}$ in $\Sigma(W)$.
\end{Example}

\section{Examples of Quiver Presentations.}
\label{sec:examples}

In this section we look at particular examples of irreducible finite Coxeter
groups.  For each of the series $A_n$, $B_n$ and $D_n$, we list some general
properties of the quiver $\Qb = (\Vb, \Eb)$ and give a presentation as a
quiver with relations of $\Sigma(W)$ for the smallest group $W$ in the series
for which the descent algebra is not a path algebra.

\subsection{Algorithm}
\label{alg:main}
Based on Theorem~\ref{thm:main} and the results
of the previous section, we can use the following algorithm to calculate
a quiver presentation for the descent algebra $\Sigma(W)$ of a
particular finite Coxeter group~$W$.

  \begin{itemize}
  \item \textbf{Given:} a finite Coxeter group $W$.
  \item \textbf{Compute:} A quiver $\Qb = (\Vb, \Eb)$ and a set $\RR$ of
    relations between the paths in $\Qb$ such that the path algebra of $\Qb$
    modulo $\RR$ is isomorphic to $\Sigma(W)$.
  \end{itemize}
  \begin{enumerate} \renewcommand{\labelenumi}{\arabic{enumi}.}
  \item $\Vb \gets \Lambda$, the set of all shapes of $W$.
  \item $M \gets \{\alpha \in \Psi: l(\alpha) > 0 \text{ and } \Delta(\alpha) \neq 0\}$;
  \item $i \gets 0$.
  \item while $M \neq \emptyset$:
  \item \qquad $i \gets i + 1$;
  \item \qquad $E_i \gets M$;
  \item \qquad add to $\RR$ a basis of the nullspace of $\Delta$ on $\Span{E_1
    \cup \dots \cup E_i}_{\Q}$;
  \item \qquad remove redundant elements from $E_1$;
  \item \qquad $M \gets M \circ E_1$;
  \item Return $(\Vb, E_1)$, $\RR$ expressed in terms of $E_1$.
  \end{enumerate}
  In the resulting quiver, the  edges are  elements
  of $\Psi$, so that an explicit  isomorphism between the path algebra of this
  quiver  $\Qb$   and  the  descent   algebra  $\Sigma(W)$  is   obtained  by
  simply applying~$\Delta$.

\subsection{Type $A_n$.}
The shapes of a Coxeter group of type $A_n$ correspond to the partitions of
$n+1$.  If $\lambda, \lambda' \in \Lambda$ are such that $\lambda' \lessdot
\lambda$ then the partition $q$ corresponding to $\lambda$ is obtained from
the partition $p$ corresponding to $\lambda'$ by joining two parts of~$p$.
And if there is an edge from $\lambda'$ to $\lambda$ in $\Qb$ then the two
parts are distinct, by Theorem~\ref{thm:edge1}.
It turns out that there is in fact an edge in $\Qb$ whenever the two parts
are distinct.  It furthermore turns out that $f_{\alpha} \in \Rad^2
\Sigma(W)$ for all $\alpha \in \Psi$ with $\ell(\alpha) > 1$.

Hence the vertices of quiver $\Qb$ correspond to the partitions of $n+1$ with
an arrow $p \to q$ between partitions $p,q$ of $n+1$ if and only if $q$ is
obtained from $p$ by joining two distinct parts of $p$.  This quiver has an
isolated vertex $1^{n+1}$ and a further isolated vertex $2^{(n+1)/2}$ if $n$
is odd.  The remaining vertices form one connected component. The descent
algebra $\Sigma(W)$ therefore has $2$ or $3$ blocks, depending on whether $n$
is odd or even.  This description of $\Qb$ has been given by Garsia and
Reutenauer~\cite{GarReut1989}.  It also follows from the results of
Blessenohl and Laue~\cite{BlessenohlLaue1996,BlessenohlLaue2002}, as pointed
out by Schocker~\cite{Schocker2004}.  A complete proof of this description in
the present framework together with a description of the relations in a
quiver presentation for $\Sigma(W)$ of type $A_n$ will be the subject of a
subsequent article.

\begin{Example}
  Consider the Coxeter group $W$ of type $A_5$ with Coxeter diagram:
  \begin{align*}
      1 - 2 - 3 - 4 - 5
  \end{align*}
  Here, and similarly in the following examples, we identify the elements of
  the set $S = \{1,2,\ldots,5\}$ with the simple reflections of $W$.  
\begin{table}[htb]
  \centering
$\begin{array}{|ccc|ccc|ccc|}
\hline
\vb & \text{type} & \lambda &
\vb & \text{type} & \lambda &
\vb & \text{type} & \lambda \\
\hline
1. & 111111 & [\emptyset] & 
5. & 222 & [135] &
9. & 42 & [1235] \\
2. & 21111 & [1] &
6. & 321 & [124] &
10. & 51 & [1234] \\
3. & 2211 & [13] &
7. & 411 & [123] &
11. & 6 & [S] \\ 
4. & 3111 & [12] &
8. & 33 & [1245] &
&& \\ \hline
\end{array}$

\medskip
$\begin{array}{|cc|cc|cc|}
\hline
\eb & \alpha &
\eb& \alpha &
\eb & \alpha \\
\hline
2 \to 4. & [12;1] &
6 \to 8. & [1245;1] &
7 \to 10. & [1234;1] \\
3 \to 6. & [124;1] &
6 \to 9. & [1235;1] &
9 \to 11. & [S;2] \\
4 \to 7. & [123;1] &
6 \to 10. & [1234;2] &
10 \to 11. & [S;1] \\
\hline
\end{array}$
  \caption{\strut The quiver of $\Sigma(W)$ for $W$ of type $A_5$.}
  \label{tab:a5}
\end{table}
The
  vertices $\vb \in \Vb$ of the quiver $\Qb$, correspond to the partitions
  of $6$ and are enumerated in  Table~\ref{tab:a5}, together with a
  representative $L \subseteq S$ for each shape $\lambda = [L]$ of $W$.
The edges $\eb \in \Eb$ of the quiver $\Qb$ are listed in terms of the 
vertex numbering and as a streets $\alpha$.
The only relation in this case is 
\begin{align}
(3
\to
6
\to
9
\to
11)
 =
(3
\to
6
\to
10
\to
11),
\end{align}
arising from the fact that $f_{[S;234]} = f_{[S;134]}$ in $\Sigma(W)$.
\end{Example}

\subsection{Type $B_n$.}
In the Coxeter group $W$ of type $B_n$ the longest element $w_0$ is central.
Proposition~\ref{pro:fLs=0} thus yields $f_{\alpha} = 0$ for
all $\alpha \in \Psi$ with $\ell(\alpha) = 1$, and therefore no cover
relation $\lambda' \lessdot \lambda$ of shapes $\lambda, \lambda' \in
\Lambda$ gives rise to an edge of the quiver $\Qb$!

\begin{table}[htb]
  \centering
$\begin{array}{|ccc|ccc|ccc|}
\hline
\vb & \text{type} & \lambda &
\vb & \text{type} & \lambda &
\vb & \text{type} & \lambda \\
\hline
1.  & 111111 & [\emptyset] &
11. & 321 & [235] &
21. & 21 & [1236] \\
2.  & 11111 & [1] &
12. & 211 & [125] &
22. & 51 & [2345] \\
3.  & 21111 & [2] &
13. & 411 & [345] &
23. & 11 & [1234] \\
4.  & 2111 & [16] &
14. & 111 & [123] &
24. & 4 & [12456] \\
5.  & 2211 & [26] &
15. & 32 & [1356] &
25. & 5 & [13456] \\
6.  & 3111 & [45] &
16. & 22 & [1246] &
26. & 3 & [12356] \\
7.  & 1111 & [12] &
17. & 33 & [2356] &
27. & 6 & [23456] \\
8.  & 222 & [246] &
18. & 31 & [1256] &
28. & 2 & [12346] \\
9.  & 221 & [146] &
19. & 41 & [1345] &
29. & 1 & [12345] \\
10. & 311 & [145] &
20. & 42 & [2456] &
30. & \emptyset & [S] \\
\hline
\end{array}$

\medskip
$\begin{array}{|cc|cc|cc|}
\hline
\eb & \alpha &
\eb & \alpha &
\eb & \alpha \\
\hline
3 \to 13. & [234;23] & 
9 \to 28. & [12346;23] &   
12 \to 29. & [12345;34] \\
3 \to 14. & [123;12] & 
10 \to 25. & [13456;34] &
13 \to 27. & [23456;23] \\  
4 \to 19. & [1345;34] &
10 \to 29. & [12345;23] &                               
13 \to 29. & [12345;12] \\   
4 \to 23. & [1234;23] & 
11 \to 26. & [12356;12] &                               
15 \to 30. & [S;24] \\ 
5 \to 20. & [2346;23] & 
11 \stackrel{.}{\to} 27. & [23456;34] &                             
18 \to 30. & [S;34] \\ 
5 \to 21. & [1235;12] &
11 \stackrel{..}{\to} 27. & [23456;42] &                               
19 \to 30. & [S;23] \\  
5 \to 22. & [2345;34] & 
11 \to 28. & [12346;12] &                               
20 \to 30. & [S;13] \\ 
6 \to 22. & [2345;23] &
11 \to 29. & [12345;13] &                             
21 \to 30. & [S;45] \\ 
6 \to 23. & [1234;12] &
12 \to 24. & [12456;45] &                               
22 \to 30. & [S;12] \\  
9 \to 25. & [13456;45] &
& & & \\ \hline
\end{array}$
  \caption{\strut The quiver of $\Sigma(W)$ for $W$ of type $B_6$.}
  \label{tab:b6}
\end{table}

The shapes of $W$ correspond to the partitions of $m \in \{0, \dots, n\}$.
Experimental evidence suggests that the edge set $\Eb$ of the quiver $\Qb$ on
this vertex set can be described as follows.  There is an $e$-fold edge $p
\to q$ between two partitions $p,q$ if $q$ is be obtained from
$p$ by either joining $3$ parts $a, b, c$ with $e = \Size{\{a, b, c\}} - 1$,
or by dropping $2$ parts $a, b$ with $e = \Size{\{a, b\}} - 1$.

The graph described by these rules has five isolated vertices $1^n$,
$1^{n-1}$, $1^{n-2}$, $2^{\Floor{n/2}}$, and $3^{\Floor{n/3}}$, for $n$ large
enough, and a further isolated vertex $2^{n/2-1}$ if $n$ is even.  The
remaining vertices form two connected components, one on the partitions of
odd length and one on the partitions of even length.  The descent algebra
$\Sigma(W)$ therefore has $4$ or $6$ blocks if $n = 2$ or $n = 3$, and, for
larger $n$, it has $8$ or $7$ blocks, depending on whether $n$ is even or
odd.  The quiver  is illustrated with an example  of type $B_6$ below.

We hope to give a complete proof of this description of the quiver together
with a description of the relations in a quiver presentation for
$\Sigma(W)$ of type $B_n$ in a subsequent article.

\begin{Example}
  Consider the Coxeter group $W$ of type $B_6$ with Coxeter diagram:
\begin{align*}
    1 = 2 - 3 - 4 - 5 - 6
\end{align*}
The shapes of $W$, which serve  as vertex set of the quiver $\Qb$, correspond
to the partitions of all $m \in \{0, 1, \dots, 6\}$ and are enumerated in   Table~\ref{tab:b6}, together  with a  representative $L  \subseteq S$  for each
shape $\lambda = [L]$ of~$W$.
The edges $\eb \in \Eb$ of the quiver $\Qb$ are listed in terms of the
vertex numbering and as a streets $\alpha$.
Note that there are two edges between vertices $11$ and $27$,
i.e., between the partitions $321$ and $6$ of~$6$,
which are distinguished by using the symbols $\stackrel{.}{\to}$
and $\stackrel{..}{\to}$ as arrows.

The only relation in this case is
\begin{align}
  (5 \to 20 \to 30) = (5 \to 22 \to 30),
\end{align}
arising from the fact that $f_{[2346;23]} f_{[S;13]} = f_{[2345;34]} f_{[S;12]}$
in~$\Sigma(W)$.
\end{Example}

\subsection{Type $D_n$.}
If $n$ is even then the longest element $w_0$ is central in the Coxeter group
$W$ of type $D_n$ and, as in the case of type $B_n$, no cover relation
$\lambda' \lessdot \lambda$ of shapes $\lambda, \lambda' \in \Lambda$ gives
rise to an edge of the quiver $\Qb$.

If $n$ is odd, then the shapes $\lambda \in \Lambda$  of $W$ correspond to the partitions of $m \in
\{0, 1,  \dots, n-2\} \cup  \{n\}$
in such a way that each part $a$ of a partition $p$ of $m$ stands for 
a direct  factor of type $A_{a-1}$ of $W_L$, $L \in \lambda$, and 
if $m < n$ then $W_L$ also has a direct factor of type $D_{n-m}$.
And if $\lambda', \lambda \in \Lambda$ correspond to 
partitions $p, q$ in this set then
$\lambda' \lessdot \lambda$ if and only if
$q$ is obtained from $p$ by either joining two parts or by dropping one part.
With Theorem~\ref{thm:edge1}, it  can be
shown  that, if 
$\lambda' \lessdot \lambda$ and there is an edge from $\lambda'$ to $\lambda$
in the quiver $\Qb$, then $\lambda'$ corresponds to a partition $p$ of $n$,
which has exactly one odd part $a$, and the partition $q$ corresponding to
$\lambda$ arises from $p$ by either joining $a$ and another (even) part of
$p$ or, if $a > 1$,  by dropping $a$ from~$p$.

\begin{table}[htb]
  \centering
$\begin{array}{|rcc|rcc|rcc|}
\hline
\vb & \text{type} & \lambda &
\vb & \text{type} & \lambda &
\vb & \text{type} & \lambda \\
\hline
1.  & 111111 & [\emptyset] &
10. & 111   & [123] &
19. & 11    & [1234] \\
2.  & 21111 & [1] &
11. & 411   & [134] &
20. & 4     & [12456] \\
3.  & 1111  & [12] &
12. & 22    & [1246] &
21. & 3     & [12356] \\
4.  & 2211  & [46] &
13. & 31    & [1256] &
22. & 2     & [12346] \\
5.  & 3111  & [13] &
14. & 33    & [1356] &
23. & 6^{-} & [13456] \\
6.  & 222^{-} & [146] &
15. & 21    & [1236] &
24. & 6^{+} & [23456] \\
7.  & 222^{+} & [246] &
16. & 42^{-} & [1346] &
25. & 1     & [12345] \\
8.  & 211   & [124] &
17. & 42^{+} & [2346] &
26. & \emptyset & [S] \\
9.  & 321   & [236] &
18. & 51    & [1345] &
&& \\\hline
\end{array}$

\medskip
$\begin{array}{|cc|cc|cc|}
\hline
\eb & \alpha &
\eb & \alpha &
\eb & \alpha \\
\hline
2 \to 10. & [123;12] &
                           8 \to 25. & [12345;34] &  
                                                      11 \to 23. & [13456;13] \\
2 \to 11. & [134;13] &                                                          
9 \stackrel{.}{\to} 21. & [12356;12] &                            
                                                      11 \to 24. & [23456;23] \\
4 \to 14. & [1356;15] &                                                       
9 \stackrel{..}{\to} 21. & [12356;15] &                             
                                                      11 \to 25. & [12345;12] \\
4 \to 15. & [1235;12] &                                                         
9 \stackrel{.}{\to} 22. & [12346;21] &                             
                                                      13 \to 26. & [S;34] \\    
4 \to 16. & [1346;13] &                                                         
9 \stackrel{..}{\to} 22. & [12346;12] &                            
                                                      14 \to 26. & [S;41] \\
4 \to 17. & [2346;23] &                                                        
9 \stackrel{.}{\to} 23. & [13456;34] &                             
                                                      15 \to 26. & [S;45] \\    
4 \to 18. & [1345;34] &                                                         
9 \stackrel{..}{\to} 23. & [13456;41] &                             
                                                      16 \to 26. & [S;23] \\    
5 \to 18. & [1345;13] &                                                         
9 \stackrel{.}{\to} 24. & [23456;42] &                            
                                                      17 \to 26. & [S;13] \\    
5 \to 19. & [1234;12] &                                                        
9 \stackrel{..}{\to} 24. & [23456;34] &                             
18 \stackrel{.}{\to} 26. & [S;21] \\    
8 \to 20. & [12456;45] &                                                        
                           9 \to 25. & [12345;13] &                             
18 \stackrel{..}{\to} 26. & [S;12] \\    
\hline
\end{array}$
  \caption{\strut The quiver of $\Sigma(W)$ for $W$ of type $D_6$.}
  \label{tab:d6}
\end{table}
\begin{Example}
  Consider the Coxeter group $W$ of type $D_6$ with Coxeter diagram:
\begin{align*}\baselineskip0pt
  2 - \vbox{\hbox{$1$}\hbox{\,$|$}\hbox{$3$}} - 4 - 5 - 6
\end{align*}
The shapes of $W$, which serve as vertex set of the quiver $\Qb$, correspond
to the partitions of all $m \in \{0, 1, 2,3,4, 6\}$, with two copies of each
even partition of $6$, and are enumerated in Table~\ref{tab:d6}, together
with a representative $L \subseteq S$ for each shape $\lambda = [L]$ of~$W$.
The edges $\eb \in \Eb$ of the quiver $\Qb$ are listed in terms of the vertex
numbering and as a streets $\alpha$.
Note that there are two edges between vertices
$9$ and $21$, $22$, $23$, $24$ respectively, and
between vertices $18$ and $26$.

There are three relations:
\begin{align}
  (4 \to 14 \to 26) = -2 (4 \to 15 \to 26),
\end{align}
arising from  $f_{[1356;15]} f_{[S; 41]} = -2 f_{[1235;12]} f_{[S;45]}$ in $\Sigma(W)$,
\begin{align}
  (4 \to 16 \to 26)  = (4 \to 18 \stackrel{.}{\to} 26),
\end{align}
arising from  $f_{[1346;13]} f_{[S;23]} = f_{[1345;34]}
f_{[S;21]}$ in $\Sigma(W)$, and
\begin{align}
  (4 \to 17 \to 26)  = (4 \to 18 \stackrel{..}{\to} 26),
\end{align}
arising from  $f_{[2346;23]}f_{[S;13]} = f_{[1345;34]}f_{[S;12]}$ in $\Sigma(W)$.
\end{Example}

\subsection{Exceptional types.}
The quivers $\Qb$ of the descent algebras
$\Sigma(W)$ for Coxeter groups $W$ of exceptional or
non-crystallographic type have been computed 
with Algorithm~\ref{alg:main} and
are described in detail elsewhere~\cite{pfeiffer-quivefhi}.

\bibliographystyle{amsplain} 
\bibliography{descent}

\end{document}